\documentclass[12pt,reqno]{amsart}
\usepackage{amssymb}
\usepackage[all]{xy}

\def\Q {{\mathbb Q}}

\def\R {{\mathbb R}}
\def\N {{\mathbb N}}

\def\F {{\mathbb F}}
\def\Z {{\mathbb Z}}

\def\C {{\mathbb C}}
\def \III {\rm III}

\newtheorem{theorem}{Theorem}

\newtheorem{proposition}{Proposition}
\newtheorem{corollary}{Corollary}
\newtheorem{lemma}{Lemma}

\newtheorem{open question}{Open Question}
\theoremstyle{definition}
\newtheorem{definition}{Definition}
\newtheorem{remark}{Remark}
\theoremstyle{remark}

\begin{document}
\title[Lefschetz and   laminations]{Scaling  group flow and Lefschetz trace formula for 
laminated spaces with $p-$adic transversal}

\author{Eric Leichtnam}
\date{\today}

\address{Institut de Jussieu et CNRS,
Etage 7E,
175 rue du Chevaleret,
75013, Paris,  France}
\email{write ericleichtnam and then \`at math.jussieu.fr}

\begin{abstract} In his approach to analytic number theory C. Deninger 
has suggested that to the Riemann zeta function $\widehat{\zeta}(s)$ (resp.  the zeta function 
$\zeta_Y(s)$ of a smooth projective curve $Y$ over a finite field 
$\F_q$, $q=p^f$))
one could possibly associate a foliated Riemannian laminated space 
 $( S_\Q , \mathcal{F}, g, \phi^t)$ 
(resp. $( S_Y , \mathcal{F}, g, \phi^t)$) endowed with an action of  a flow $\phi^t$ 
whose primitive compact orbits should correspond to the primes of 
$\Q$ (resp. $Y$). Precise  conjectures were stated in our report 
\cite{El} on Deninger's work.
The existence of  such  a foliated space and flow $\phi^t$ 
is still unknown except when $Y$ is an elliptic curve (see Deninger \cite{D 2}). 
Being
motivated by this latter case, we introduce a class of foliated laminated spaces
$( S=\frac{\mathcal{L}\times \R^{+*}}{q^\Z} , \mathcal{F}, g, \phi^t)$ where $\mathcal{L} $ is locally $D\times \Z_p^m$, $D$ being an open disk of $\C.$  Assuming 
that the leafwise harmonic forms on $ \mathcal{L}$ are 
locally constant transversally, we prove a Lefschetz trace formula 
for the flow $\phi^t$ acting on the leafwise Hodge cohomology $H^j_\tau$ ($0\leq j \leq 2$)
of $(S,\mathcal{F})$ that is very similar to the explicit formula 
for the zeta function of a (general) smooth curve over $\F_q$. 
We also prove that the eigenvalues of the infinitesimal 
generator of the action of $\phi^t$ on $H^1_\tau$ have real 
part equal to $\frac{1}{2}.$

Moreover,  we suggest in a precise way  that the flow $\phi^t$ should be 
induced by a renormalization group flow "\`a la K.  Wilson". We show that when 
$Y$ is an elliptic curve over $\F_q$ this is indeed the case. 
It would be very interesting to establish 
a precise connection between  our results and those of   
Connes (page 553 \cite{C3}, page 90 \cite{C4}) and Connes-Marcolli \cite{CM1}, \cite{CM2} on the Galois interpretation of the 
renormalization group.

\end{abstract}
\maketitle

\tableofcontents
\section{Introduction}$\;$

Several papers (see \cite{D 1}, \cite{D 0c}, \cite{D 2}, \cite{D 1c}, \cite{D 2b}) of 
Deninger  lead to suggest that to the Riemann zeta function 
$\widehat{\zeta}(s)$ one could possibly associate 
the following two data:

1) A Riemannian foliated space 
of the form
$$
( \overline{S_\Q}=\frac{ \mathcal{L} \times \R^{+*}}{\Q^{+*}} \cup \frac{ \mathcal{L}}{\Q^{+*}}, \mathcal{F},g)
$$ where $\mathcal{L}$ is a $\sigma-$compact complex $1-$dimensional laminated space 
on which $\Q^{+*}$ acts. The path connected components of $\mathcal{L}$ induce 
a foliation of  $\frac{ \mathcal{L} \times \R^{+*}}{\Q^{+*}} $ by Riemann surfaces and $g$ is a leafwise kaehler
metric.

2) A flow $\phi^t$ acting on $( \overline{S_\Q}=\frac{ \mathcal{L} \times \R^{+*}}{\Q^{+*}} \cup \frac{ \mathcal{L}}{\Q^{+*}}, \mathcal{F} )$ whose primitive closed orbits correspond to the primes 
of $\Q$ and admitting a fixed point in $ \frac{ \mathcal{L}}{\Q^{+*}}$.
The action of $\phi^t$  on the $\mathcal{L}-$leaf space 
$\frac{ \R^{+*}}{\Q^{+*}} \cup \{ {\rm pt}\}$ should be given 
by $\phi^t ([x]) = [ e^{-t} x ]$ and $\phi^t ({\rm pt} ) ={\rm pt}.$ Moreover $(\phi^t)^*[\lambda_g]=e^t [\lambda_g]$ where $[\lambda_g]$ denotes  the reduced leafwise cohomology class of the 
leafwise kaehler form associated to $g.$

\medskip
The quotient $ \frac{ \mathcal{L}}{\Q^{+*}}$ allows to compactify the space $ \frac{ \mathcal{L} \times \R^{+*}}{\Q^{+*}}$.  For more precise informations and axioms see \cite{El} or even Section \ref{Section:curve}. 
Notice that the existence of such a quadruple $( \overline{S_\Q}, \mathcal{F}, g, \phi^t )$ is still unknown.  Recall that Alain Connes \cite{C 2}  has reduced the validity of the Riemann hypothesis 
to a trace formula (of Lefschetz type) on the quotient 
space $ \frac {\mathbb{A}_f \times \R}{\Q^*}.$ Notice nevertheless that 
$ \frac {\mathbb{A}_f \times \R}{\Q^*}$ does not satisfy the properties required for 
$( \overline{S_\Q}, \mathcal{F}, g, \phi^t ).$

It is  natural to guess that the action of $\phi^t$
on $\overline{S_\Q}$ should be induced by the action on 
$\mathcal{L}\times \R^{+*}$ of a flow still denoted $\phi^t$. 
Given the previous property 2) above,  one is lead to 
search a flow of the form
$$
 \phi^t(l,x)= (\psi^t_x(l), x e^{-t} ),\; \forall (l,x) \in \mathcal{L}\times \R^{+*}.
$$ But if we write $\psi^t_x= R_{x , x \,e^{-t}}$ 
then we recognize the general scheme of the method of renormalization group  " \`a la K. Wilson" ([page 554]\cite{Del}),  where $R_{x , x \,e^{-t}} $ should act 
on a space $\mathcal{L}$ of lagrangians.
We noticed this point in 2002 when 
we were preparing our lectures (on which \cite{El} is based) 
on Deninger's work  for Bar Ilan University.

Then we are lead to guess that  
a similar picture should exist for the zeta function of a smooth projective curve $Y$ over  
$\F_q.$ Namely that there should exists a Riemannian foliated laminated space 
$(S_Y, \mathcal{F}, g, \phi^t)$ where $\phi^t$ should be a renormalization group flow
whose primitive closed orbits should correspond to the closed points of $Y.$

In Section 2 we consider the case of an elliptic curve $E_0$ and the 
 Riemannian foliated laminated space 
$(S(E_0), \mathcal{F}, g, \phi^t)$ constructed by Deninger. Then we show that 
$(S(E_0), \mathcal{F}, g, \phi^t)$ can be interpreted as a renormalization group flow.  
Now we motivate our main result and explain its main context. 
Recall that the following explicit formula is satisfied for 
the zeta function $\zeta_Y(s)$ of $Y.$ Let $\alpha \in C^\infty_c(\R,\R).$ Then:

\begin{equation} \label{eq:exp} \begin{aligned}[rcl]  
\sum_{\nu \in \Z} 
& \int_\R \alpha(t) e^{t\, ( \frac {   2 i \pi \nu }{\log q}) }  d t -
\sum_{j=1}^{2 g}\sum_{\nu \in \Z}   \int_\R \alpha(t) e^{t\, ( \rho_j+\frac {   2 i \pi \nu }{\log q}) }  d t +
 \sum_{\nu \in \Z}   \int_\R \alpha(t) e^{t\, ( 1 + \frac {  2 i \pi \nu }{\log q}) }  d t =
\\ &
(2-2 g)\alpha(0) \log q + \sum_{\gamma} \sum_{k \geq 1} l(\gamma)
\left( \,e^{- k l(\gamma) } \alpha( -k l(\gamma) ) + \alpha( k l(\gamma) )\, \right)
\end{aligned}
\end{equation}
 where the $(\rho_j+\frac {   2 i \pi \nu }{\log q} )$ run over the zeroes 
of $\zeta_Y$ and the $\gamma$ (with norm $e^{l(\gamma)}$) run over the closed points 
of $Y$. As pointed out by Deninger, the dissymmetry of the coefficients of 
$\alpha(-k l(\gamma))$ and $ \alpha(k l(\gamma))$ in the right handside is of arithmetic narure. It arises when one tries to intertwine the functional equation (addition) 
and the Eulerian product (multiplication) in the proof of the  explicit formula \eqref{eq:exp}.

Deninger has suggested that it might be possible to interpret 
the above formula as a Lefschetz trace formula. Having in mind 
in fact the case of number fields,   he made interesting remarks on dynamical Lefschetz trace formulas on 
laminated foliated spaces see [Section5]\cite{D 2b}.

Consider briefly the case of a compact connected three dimensional manifold $X$ endowed with
a codimension one  foliation $(X,\mathcal{F})$. Assume that $(X,\mathcal{F})$ is endowed 
with a flow $\phi^t$ which acts transversally and whose  
 closed orbits of $\phi^t$ are simple, thus  $(X,\mathcal{F})$ is Riemannian.  Let $\pi^j_\tau$ denote the projection onto the 
leafwise Hodge cohomology $H^j_\tau$ $(0\leq j \leq 2$). Then  Alvarez-Lopez and 
Kordyukov \cite{A-K2} have proved  the following Lefschetz trace formula:
$$
\sum_{j=0}^2 (-1)^j {\rm TR} \int_\R \alpha(s) \,\pi^j_\tau\circ (\phi^s)^* \circ \pi^j_\tau \,d s\,=
$$
$$\chi_\Lambda \alpha(0) + \sum_{\gamma} \sum_{k \geq 1} l(\gamma)
\left( \, \epsilon_{-k} \alpha( -k l(\gamma) ) + \epsilon_k \alpha( k l(\gamma) )\, \right)
$$ where 
$\gamma$ runs over the primitive closed orbits of $\phi^t$,  
$\epsilon_{\pm k}= {\rm sign} \,{\rm det} ({\rm id}- D\phi_{| {\rm T}_y \mathcal{F}}^{\pm k l(\gamma)}),$ $y\in \gamma$ and ${\rm TR}$ denotes the (usual) trace of  trace class
operators. Notice  that here there is no dissymmetry for the coefficients 
of $\alpha( -k l(\gamma) )$ and $\alpha( k l(\gamma) ).$
The reason for this absence of dissymmetry is due to the Guillemin-Sternberg formula \cite{G-S} which states that the geometric contribution of a closed orbit $\pm k \gamma$ should be:
$$
l(\gamma) \alpha(\pm k l(\gamma) )\, \sum_{j=0}^2 (-1)^j\,
\frac{ {\rm Tr} \bigl( (D \phi^{\pm k l(\gamma)})^*: \wedge^j {\rm T}^*_y  \mathcal{F} \mapsto \wedge^j {\rm T}^*_y  \mathcal{F} \bigr)}
 {|{\rm det}( {\rm id } - D \phi^{\pm k l(\gamma)}_{| {\rm T}_y \mathcal{F}}) |}
$$ where $y$ is any point on $\gamma.$
In Lemma \ref{lem:trans} of Section 2.4 we shall provide 
in the case of an elliptic curve $E_0$ over $\F_q$  a dynamical 
explanation of this dissymmetry pointing out the role 
of the $p-$adic transversal in $S(E_0)$. This will lead us to propose a list 
of four Assumptions in Section \ref{Section:curve} for a Riemannian foliated laminated space 
$$
(S=\frac{\mathcal{L}\times \R^{+*}}{q^\Z}, \mathcal{F}, g, \phi^t)
$$ where $\mathcal{L}$ is locally of the form 
$D\times \Z_p^m$, $D$ is a disk of $\C.$ The leaves of $S$ are induced 
by $\mathcal{C}\times \{x\}$ where $\mathcal{C} $ is a path connected component of 
$\mathcal{L}$.
 Assumption (iv) states that the elements of the  vector space $\mathcal{H}^1_{\mathcal{L}}$  of continuous and leafwise harmonic forms 
on $\mathcal{L}$  are locally constant along the $p-$adic transversal $\Z_p^m.$ 
It implies that the vector space $\mathcal{H}^1_{\mathcal{L}}$  is of finite even dimension $2 g$. This Assumption is of course satisfied by the foliated space $(S(E_0),  \mathcal{F})$ of Section 2 but we do not know if there are any 
example which satisfies this Assumption with $g \geq 2$, see the remark  before Lemma \ref{lem:transv}.  Anyway, we are forced 
to use this Assumption for the following two reasons. 
First if we want for $(S=\frac{\mathcal{L}\times \R^{+*}}{q^\Z}, \mathcal{F}, g, \phi^t)$ a Lefschetz trace formula of the type \eqref{eq:exp} then it seems
that we need $\mathcal{H}^1_{\mathcal{L}}$  to be of finite dimension $2 g.$
Second, we shall need the following operator 
$$
\int_\R \alpha(s) (\phi^t)^* d s \circ \pi^j_\tau
$$ to be trace class. Morally to $\pi^j_\tau$ (resp. $ \int_\R \alpha(s) (\phi^t)^* d s$) corresponds a regularizing process along the leaves $[\mathcal{C}\times \{x\}] $ of $S$ (resp. the 
integral curves of $\phi^t$). But $ \int_\R \alpha(s) (\phi^t)^* d s \circ \pi^j_\tau$ does not involve any regularizing process along 
the $p-$adic transversal $\Z_p^m$. Therefore,  unless Assumption (iv) is satisfied we see no reason 
why $ \int_\R \alpha(s) (\phi^t)^* d s \circ \pi^j_\tau$ should be trace class.

Then, our main result is Theorem \ref{Thm:trace}.  It  proves a Lefschetz trace formula similar to \eqref{eq:exp} and shows that the eigenvalues of the infinitesimal generator of $\pi^1_\tau\circ (\phi^t)^*$ acting on the leafwise Hodge cohomology $H^1_\tau$ 
have real part equal to $\frac{1}{2}.$  Our main new ingredient is the transversal $p-$adic 
Laplacian $\Delta_{p,T}$ on $\mathcal{L}$ (see Definition \ref{def:Delta}). The intrinsic meaning of Assumption (iv) is  the inclusion $\mathcal{H}^1_{\mathcal{L}} \subset \ker \Delta_{p,T}.$
Moreover, Assumption (iv) allows 
to use a "contraction process" along the $p-$adic transversal of $S$ (See Definition \ref{Def:contraction}). 
Then we use in essential way 
results of Alvarez-Lopez and Kordyukov since in some sense our 
foliated laminated spaces are closed to Riemannian foliations.
Even in the case of an elliptic curve 
$E_0$, it seems of some interest to provide a proof of 
\eqref{eq:exp} \`a la Atiyah-Bott-Lefschetz and the new ingredients introduced here should be useful in other 
contexts.

We also show (Proposition \ref{Wilson}) that the von Neumann algebra 
$W(S,\mathcal{F})\rtimes_{\phi^t} \R$, which describes the non commutative space 
of closed points, is of type ${\rm III}_{\frac{1}{q}}.$ This matches with 
Connes's approach (see \cite{C 2}, \cite{El}).

In a future paper \cite{El1} we shall try to  propose a list of axioms that
a Riemannian foliated laminated space
$( \overline{S_\Q}, \mathcal{F}, g, \phi^t )$ should satisfy in order to get 
a Lefschetz trace formula that should be analogous to the explicit formula
for the Riemann Zeta function.

\medskip Now we try to explain why there should exist a  connection with the work 
of Connes-Marcolli.
Let $K$ be a local field with residue class field  the finite field $\F_q,$ 
let $f: z \rightarrow z^q$ be the canonical generator of 
the Galois group of $\overline{\F_q}$ over $\F_q$. Then 
local class field theory shows that the Galois group of the maximal unramified 
extension of $K$ admits a dense subgroup which is naturally isomorphic 
to  $\{ f^k,\;  k \in \Z \}.$ Now let $Y$ be a smooth projective variety over $\F_q$, 
then (see [page 292]\cite{Milne})  the automorphism 
${\rm Id }\otimes f$ of $\overline{Y}=Y \times_{\F_q} \overline{\F_q}$ is called   the 
arithmetic Frobenius. On the other hand, 
the geometric Frobenius $F: \overline{Y} \rightarrow  \overline{Y}$ is the 
$\overline{\F_q}-$morphism sending the point $P$ with coordinates 
$(a_i), \, a_i \in \overline{\F_q},$ to the point $F(P)$ with coordinates 
$(a_i^q).$
The action of $({\rm Id} \otimes f)^{-1}$ on $l-$adic cohomology coincides 
with the action of the geometric Frobenius $F$ (see [page 292]\cite{Milne}).

On the other hand, A. Connes 
(page 553 \cite{C3}, page 90 \cite{C4}) has suggested that $\R^{+*}$, as part of the renormalization group, 
should play the role of the missing (unramified) Galois group at the archimedean place of $\Q$. 
A. Connes is motivated by an unramified local class field analogy and his classification of type 
$\III-$factors. Moreover, Connes-Marcolli  \cite{CM2} have shown that 
the renormalization group flow is an ambiguity Galois group acting 
on Quantum Field Theories (QFT's). This   seems  reminiscent of 
a continuous version of ${\rm Id} \otimes f$.

 We suggest that the scaling group flow $\phi^t$, in 
$(S,\mathcal{F},g, \phi^t)$ as above, 
corresponds to a continuous version of $F$ and that
Connes's suggestion  corresponds to a continuous version 
of ${\rm Id} \otimes f.$ 
 It would be interesting to establish a precise connection between these 
two continuous notions of Frobenius by providing an ambiguity galois group 
interpretation of $\phi^t.$ 
The  results of \cite{CM1}, \cite{CM2} should be very helpful with this respect. 
Moreover, this should 
 allow to decide if to a smooth curve $Y$ over $\F_q$ one can (or cannot)
associate a foliated space $(S_Y, \mathcal{F}, g, \phi^t)$ satisfying the four 
Assumptions of Section \ref{Section:curve} and such that the primitive closed orbits 
of $\phi^t$ should correspond to the closed points of $Y$. 
It is also tempting to try to establish a connection with Haran's recent approach \cite{Haran}.

\section{The case of an elliptic curve $E_0$ over $\F_q$} $\;$

\subsection{The zeta function $\zeta_{E_0}(s)$ and the explicit formula}$\;$

\medskip
Let $E_0$ be an elliptic curve over a finite field $\F_q$. Recall that
the zeta function $\zeta_{E_0}(s)$ of $E_0$ is given by:
\begin{equation}  \label{zetaE}
\zeta_{E_0}(s)= \prod_{w \in |E_0|} \frac {1} {1-(Nw)^{-s}} \,=\,
\frac {(1-\xi q^{-s})(1-\overline{\xi }q^{-s})} 
{(1-q^{-s})(1-q^{1-s})}
\end{equation}
 where $|E_0| $ denotes the set of closed points of 
$E_0$  and  $\xi$ is a complex number which by Hasse's theorem satisfies
$|\xi|= \sqrt{q}$. The explicit formula for $ \zeta_{E_0}(s)$ takes
the following  form. Let $\alpha \in C^\infty_c(\R, \R)$ and set for
any real $s$, $\Phi(s)= \int_\R e^{s t} \alpha(t)\,dt$. Then, one has:
\begin{equation} \label{EF2} 
\begin{aligned} [rcl]
\sum_{\nu \in \Z} 
&
\Phi
\left( \frac {2 \pi \nu i} {\log q } \right)  
-
\sum_{\rho \in \zeta_{E_0}^{-1}\{ 0 \}} \Phi (\rho)
+\sum_{\nu \in \Z} \Phi \left( 1 + \frac {2 \pi \nu i} {\log q } \right) = 
\\ &
=
\sum_{w \in | E_0 |} \log N w \left( 
\sum_{k \geq 1} \alpha (k \log N w) + 
\sum_{k \leq  -1} (N w)^k  \alpha (k \log N w) \right) .
\end{aligned}
\end{equation}
The idea of the proof is to apply the residue theorem to 
$$
s\rightarrow \left(\int_0^{+\infty} \sqrt{t} \,\alpha(\log t) t^s \frac{d\,t}{t}\right)\, \frac{\zeta_{E_0}^\prime}{\zeta_{E_0}}(s)
$$ and to use the functional equation $\zeta_{E_0}(s)=\zeta_{E_0}(1-s)$.
At the end of this Section we shall explain briefly how Deninger managed (see \cite{D 2} for the 
details) to interpret this explicit formula \eqref{EF2} as a Lefschetz trace 
formula.

\subsection{The Riemannian laminated foliated space $(S(E_0), \mathcal{F}, g, \phi^t)$ }$\;$

$\medskip$

Let $\phi_0: E_0 \rightarrow  E_0$ be the $q-$th power Frobenius endomorphism of $E_0$ over $\F_q$.
Deninger has used (see \cite{D 2}) the following result due to Oort \cite{Oort}:

\begin{lemma}  \label{lemmalift} There exists:
{\item {1]}}  a complete local integral domain $R$ with field of fractions $L$ a finite 
extension of $\Q_p$ ($q=p^r$) such that $R/{\mathcal{M}}=\F_q$ where ${\mathcal{M}}$ is the maximal 
ideal of $R$.
{\item  {2]}} an elliptic curve ${\mathcal{E}}$ over spec $R$  together with an endomorphism 
$\phi: {\mathcal{E}}\rightarrow {\mathcal{E}}$ such that:
$$
({\mathcal{E}}, \phi) \otimes \F_q= (E_0, \phi_0).
$$ So $ ({\mathcal{E}}, \phi)$ is a lift of $(E_0, \phi_0)$ in characteristic zero.
\end{lemma}

\begin{remark}  {\item1)} If the elliptic curve $E_0$ is ordinary, then one may take for $R$ the ring of 
Witt vectors of $\F_q$, $W(\F_q)$, and then there is a canonical choice of the lifting 
 $( {\mathcal{E}}, \phi)$. On the contrary, if $E_0$ is supersingular ([page 137]\cite{Silverman}), then there is no canonical 
choice of $( {\mathcal{E}}, \phi)$.
{\item 2)} It is  possible to lift  a curve of genus $\geq 2$ (over $\F_q$) in characteristic zero, but 
Hurwitz's formula ([page 41] \cite{Silverman}) shows that one cannot lift 
its Frobenius morphism.
\end{remark}

Now (still following  \cite{D 2}), we denote by $E={\mathcal{E}}\otimes_R L$ the generic fibre. Then 
${\rm End}_L(E)\otimes \Q=K$ is a field $K$ which is either $\Q$ or an imaginary quadratic extension 
of $\Q$. We fix an embedding $L\subset \C$ and consider the complex analytic elliptic curve 
$E(\C)$. Let $\omega $ be a non zero holomorphic one form on $E(\C)$ and let $\Gamma$ be its 
period lattice. Then the Abel-Jacobi map:
$$
E(\C) \rightarrow \C/\Gamma,\; p\rightarrow \int_0^p \omega \, {\rm mod}\, \Gamma
$$ induces an isomorphism. Next  we choose the embedding $K\subset \C$ such  that 
for   any $\alpha \in K$, $\Theta(\alpha)$ induces the 
multiplication by $\alpha$ on the Lie algebra $\C$ of $\C / \Gamma$ where $\Theta$ is 
the 
natural homomorphism:
$$
\Theta:\, K= {\rm End}_L(E)\otimes \Q \rightarrow  {\rm End} (\C / \Gamma)\otimes \Q.
$$ Next  we consider the unique element $\xi \in \Theta^{-1}({\rm End}_L(E) )\subset K$ such that 
$\Theta(\xi)= \phi \otimes L$. By construction one has 
 $\xi \Gamma \subset \Gamma$ and the complex
elliptic curve $\C/\Gamma$ endowed with the multiplication by $\xi$
represents a lift of $(E_0, \phi_0).$  Now, we set 
\[
V= \cup_{n\in \N} \xi^{-n} \Gamma, \;\; {\rm T} \Gamma =
\lim_{+\infty\leftarrow n} \frac {\Gamma} {\xi^n \Gamma},\;\;
\hbox{and}\;\; V_\xi \Gamma = {\rm T}\Gamma \otimes_\Z \Q.
\] The set ${\rm T} \Gamma$ is a Tate module defined by a projective limit and 
$V_\xi \Gamma$ is a $\Q_p-$vector space of dimension $1$ (resp. $2$) 
if $\xi \notin \Z$ (resp $\xi \in \Z$).

Any element $v$ of $V$ acts on $\C \times V_\xi \Gamma $ by $v.(z,
\hat{v})= (z+v, \hat{v}-v)$, we denote by $\frac {\C \times V_\xi
\Gamma} {V}$ the quotient space. 
\begin{lemma} \label{equiv}
{\item 1)} Let $F$ be a finite set of representatives in $\Gamma$ of the quotient 
group $\frac{\Gamma}{\xi \Gamma}.$ Then any element of 
$V_\xi \Gamma$ is of the form 
$\sum_{l\geq -k} a_l \xi^l$ where $k \in \N$ and  the $a_l \in F.$ Moreover the 
multiplication by $\xi$ defines an automorphism of $V_\xi \Gamma.$

{\item 2)} The natural homomorphism:
$$
\frac{\C \times {\rm T} \Gamma}{\Gamma} \rightarrow \frac {\C \times V_\xi
\Gamma} {V}
$$ defines a $\{\xi^l,\; l\in \Z\}-$equivariant isomorphism where the action of $\xi$ is induced    by the diagonal action  
on $\C \times {\rm T} \Gamma $ and $\C \times V_\xi
\Gamma $ respectively.
\end{lemma}
\begin{proof} 1) Observe that any element of 
${\rm T} \Gamma$ is of the form $\sum_{l\in N} a_l \xi^l$ where   the $a_l \in F.$ 
Using Bezout's theorem one checks that a prime number $\widehat{p}$ not dividing $q$
(ie $\not= p$) induces by multiplication an automorphism of 
${\rm T} \Gamma.$ If $\xi \notin \Z$ then the elliptic curve $\C/\Gamma$ has 
complex multiplication and its endomorphism ring is invariant under complex 
multiplication. So, in all cases,  we  have $\overline{\xi} \Gamma \subset \Gamma.$
Then using the equality $\overline{\xi} \xi =q$ one gets the results of 1). Part 2) is now easy and left to the 
reader.
\end{proof}

Now,  any element $q^\nu \in q ^\Z$
acts on $\frac {\C \times V_\xi \Gamma} {V} \times \R^{+\ast}$ by
\[
q^\nu. ( [z, \hat{v}], x ) = ( [\xi^{\nu} z, \xi^{\nu} \hat{v}],
xq^{\nu} ).
\]
 In \cite{D 2}, 
Deninger has introduced the  (compact) laminated Riemannian foliated space
 $(S(E_0) , \mathcal{F} )$  where
$$
S(E_0)\,=\, \frac {\C \times V_\xi \Gamma} {V}
\times_{q^\Z}\R^{+\ast},
$$ and the leaves of $\mathcal{F}$ are the images of the sets 
$\C \times \{\hat{v}\} \times \{x\}$ by the natural map
$\pi: \C \times V_\xi \Gamma  \times\R^{+\ast} \rightarrow S(E_0).$ Observe 
that the domain of  a  typical foliation chart is locally isomorphic to 
$D \times \Omega \times ]1,2[$ where $D$ is an open disk of $\C$, 
$\Omega$ is an open subset of ${\rm T}\Gamma$ so that the  leaves are given 
by $D \times \{\omega \}\times \{x\}$ for $(\omega,  x) \in \Omega \times
]1,2[$; the term "laminated" refers to the fact that the local 
transversal to the foliation $\mathcal{F}$ is the disconnected space 
$\Omega \times ]1,2[$. 
\begin{remark} Using the fact that $V$ (resp. $q^\Z$) acts freely on 
$V_\xi \Gamma$ (resp. $\R^{+*}$), the reader will check that  $(S(E_0) , \mathcal{F} )$
has trivial holonomy.
\end{remark}
One defines a  flow $\phi^t$ acting on $(S(E_0) , \mathcal{F} )$ and sending 
each leaf into another leaf by:
 $\phi^t (  z, \hat{v}, x )=  (  z, \hat{v}, xe^{-t} )$. 
Let $\mu_\xi$ denote a Haar measure on the group 
$V_\xi \Gamma$ then,  one has the following
\begin{lemma} (\cite{D 2}) \label{Haar} {\item 1)}
 The measure
$$
dx_1 dx_2 \otimes \mu_\xi \otimes \frac {dx} {x}
$$ on $\C \times V_\xi \Gamma \times \R^{+\ast}$ induces a measure 
$\mu$ on $S(E_0)$.
{\item 2)} The measure $\mu$ is invariant  under the action of 
$\phi^t$.
\end{lemma}
\begin{proof} {\item 1)} We just have to check that for any $\nu \in \N^*$ and any borel subset $A$ of 
$V_\xi \Gamma$ , one has 
$$\mu_\xi(\xi^\nu A)= |\xi|^{-2\nu} \mu_\xi( A)=q^{-\nu} \mu_\xi( A).
$$ Since $(\xi^{\nu})_*\mu_\xi$ is also a  
Haar measure on $V_\xi \Gamma$ it suffices to check this equality for $A={\rm T} \Gamma$. But 
this is an immediate consequence of the fact that
$$
{\rm T} \Gamma/(\xi^\nu {\rm T} \Gamma) \simeq \Gamma/(\xi^\nu \Gamma)
$$ has $|\xi|^{2\nu}=q^\nu$ elements. 
{\item 2)} This is obvious.
\end{proof}

Using the 
 fact that $|\xi|=\sqrt{q}$, one checks that the Riemannian 
metric on the bundle ${\rm T}\C \times V_\xi \Gamma  \times\R^{+\ast} $ given by:
$$
{g}_{z,\hat{v},x} (\eta_1, \eta_2)= x^{-1} { \rm Re}\, ( \eta_1 \overline{\eta_2} )
$$ induces a Riemannian metric $g$ along the leaves of $(S(E_0) , \mathcal{F} )$ so that the following property 
is satisfied:
\begin{equation} \label{disym}
\forall \eta \in {\rm T}_{[z,\hat{v},x]} \mathcal{F},\; 
g(  {\rm T}_{[z,\hat{v},x]} \phi^t(\eta), {\rm T}_{[z,\hat{v},x]} \phi^t(\eta) )
= e^t g(\eta, \eta).
\end{equation}

\subsection{Interpretation of $(S(E_0) , \mathcal{F},g, \phi^t )$ as a renormalization group flow}
\label{subsection:interp} $\;$

\medskip
Now,
recall that  the additive group $(\C,+)$ is identified to its dual $(\widehat{\C}, \times)$ in the following way.
For each character $\chi \in \widehat{\C}$ one can find 
one and only one complex number $\alpha$ such that:
\begin{equation} \label{char}
z \in \C \,\rightarrow < \chi ; z >\,=\, e^{\sqrt{-1}( \Re \alpha \Re z+   \Im \alpha  \Im z) } \in S^1.
\end{equation}
 We set $G= \frac{ \C \times {\rm T} \Gamma} { \Gamma}$, 
since $G$ is a quotient of $  \C \times {\rm T} \Gamma$, one checks that 
$\widehat{G}$ can be  identified to a subgroup of 
$\C \times \widehat{{ \rm T} \Gamma}$. Let $\Gamma^*$ denote the 
dual lattice of $\Gamma$:
$$
\Gamma^*=
\{ z \in \C/\; \forall \gamma \in \Gamma,\; 
(z ;\gamma)=  \Re z \Re \gamma +  \Im z \Im \gamma \in 2 \pi \Z\, \}.
$$
\begin{lemma} \label{isom} {\item 1]}
One has a natural group isomorphism:
$ \widehat{{\rm T} \Gamma} \simeq \frac {\cup_{n \in \N}\, (\xi^{n} \Gamma)^*}{ \Gamma^*}.$
{\item 2]} Let $\pi$ denote the projection map 
$\pi:\cup_{n \in \N}\, (\xi^{n} \Gamma)^* \mapsto 
 \frac {\cup_{n \in \N}\, (\xi^{n} \Gamma)^*}{ \Gamma^*}.$ Then one has $\widehat{G}=\{( z, \pi(z) ) \in \cup_{n \in \N}\, (\xi^{n} \Gamma)^*  \times \widehat{{ \rm T} \Gamma}\, \}.$ 

\end{lemma}
\begin{proof} 1] We fix a set $F$ of representatives in $\Gamma$ of 
the quotient group $\frac{\Gamma}{\xi \Gamma}$ which has exactly $q=|\xi|^2$ elements. Notice 
that for any $(a_k)_{k\in \N} \in F^\N$ the  series $\sum_{k\in \N} a_k \xi^k$ converges in 
 ${\rm T} \Gamma$ and that each element of ${\rm T} \Gamma$ is of the form 
$\sum_{k\in \N, a_k \in F} a_k \xi^k.$ Recall that $\Gamma$ is dense in the compact abelian group 
$ {\rm T} \Gamma$ 
and that a fundamental system of open neighborhoods of $0\in {\rm T} \Gamma$ is provided 
by the subsets $\mathcal{O}_n= \{ \sum_{k\geq n} a_k \xi^k/\; \forall k \geq n,\; a_k \in F \}$.
Now, let $\chi$ be any character  of $\Gamma$ then there exists $\alpha \in \C$ such that 
$$
\forall \gamma \in \Gamma,\;\; < \chi ; \gamma > = e^{\sqrt{-1}(\Re \alpha \Re \gamma +   \Im \alpha  \Im \gamma) } \in S^1.
$$ In fact $\alpha$ is defined in $\frac{\C}{\Gamma^*}.$ Notice that $\chi$ extends to a character of 
${\rm T} \Gamma$ if and only if $\chi\equiv 1$  on $\mathcal{O}_n$ for a suitable $n>>1$. Therefore 
$\chi$ defines a character of ${\rm T} \Gamma$ if and only if $\alpha \in (\xi^n \Gamma)^*$ for a suitable 
$n>>1$ and the result follows.

\noindent 2] An element of $\widehat{G} $ is defined by a couple $(\chi_1, \chi_2) \in 
\widehat{\C} \times  \widehat{{\rm T} \Gamma}$ such that 
$(\chi_1)_{|\Gamma}= (\chi_2)_{|\Gamma}$ and  $\chi_1$ is defined by $\alpha \in \C$ as 
in \eqref{char}. According to part 1], $(\chi_1)_{|\Gamma}$ extends to a character of 
${\rm T} \Gamma$ if and only if $\alpha \in  \cup_{n \in \N}\, (\xi^{n} \Gamma)^*$ and  the result follows.
\end{proof}

\begin{corollary} \label{G}
 One has the following group isomorphism:
 $$
\frac{ \C \times  \widehat{{\rm T} \Gamma}} {\widehat{G}}
\simeq \frac{\C}{ \Gamma^*}.
$$
\end{corollary}
\begin{proof} We fix a set $\mathcal{R} \subset \cup_{n \in \N}\, (\xi^{n} \Gamma)^*$ of 
representatives of $ \frac {\cup_{n \in \N}\, (\xi^{n} \Gamma)^*}{ \Gamma^*} $. Let 
$(u,v) \in \C \times \widehat{{\rm T} \Gamma}$, so we can find $z\in \mathcal{R}$ such that 
$\pi(z)=-v.$ According to Lemma \ref{isom}. 2], $(z, \pi(z)) \in \widehat{G}$ and 
$(z, \pi(z))\cdot (u,v) = (u + z, 0).$ Then for any $\gamma \in \Gamma^*$ one has:
 $(\gamma, \pi(\gamma) )\cdot (u +z, 0)=(u+z+ \gamma, 0).$ Now the result 
follows immediately.

\end{proof} Now  recall that 
$$S(E_0) = \frac{\C \times {\rm T} \Gamma}{ \Gamma}\times_{q^\Z} \R^{+*}
$$ and the action of $\phi^t$ on $S(E_0)$  is given by 
$\phi^t([ z , \widehat{v}, x]) = [ z , \widehat{v}, x e^{-t}].$
The following Proposition provides an interpretation of 
$(S(E_0) , \mathcal{F},g, \phi^t )$  as a renormalization group flow for suitable
quantum field theories defined over $\frac{\C}{\Gamma^*}$
  \begin{proposition} \label{Free} The set $G= \frac { \C \times {\rm T} \Gamma}{\Gamma}$ defines in a natural way a set of free lagrangians 
on $\frac{\C}{\Gamma^*}$ (ie on which the renormalization semi-group $R_{x,x e^{-t}}$ 
acts trivially, see the Appendix).
\end{proposition}
\begin{proof} Observe that 
$\widehat{G} $ is a discrete subgroup of $\C \times  \widehat{{\rm T} \Gamma}$ and thus acts on  $\C \times  \widehat{{\rm T} \Gamma}$ by translation.  According to Corollary  \ref{G}, 
$\C \times  \widehat{{\rm T} \Gamma}$ appears as a $\widehat{G}- $principal bundle 
over the  space $ \frac{ \C \times  \widehat{{\rm T} \Gamma}} {\widehat{G}}\simeq 
\frac{\C}{  \Gamma^*}.$
 Let $h \in 
\widehat{\widehat{G}} \simeq G$, then consider  the fiber product 
$\mathcal{E}_{h}=\left( \C \times  \widehat{{\rm T} \Gamma}\right) \times_{\widehat{G}} \C$ 
where for any $ \chi_1 \in \widehat{G},$ and $( (z,\widehat{v}), u ) \in \left( \C \times  \widehat{{\rm T} \Gamma}\right) \times \C$, the point   $(\, (z,\widehat{v}), u )$ is identified with 
the point $(\chi_1 \cdot (z,\widehat{v}), \chi_1(h)^{-1} u ).$  The space $\mathcal{E}_{h}$ defines 
a complex line bundle over $\frac{\C}{\Gamma^*}$ and is endowed with a natural flat
connection.  Let  $S_\pm $ denote the spinor bundle on $\frac{\C}{\Gamma^*}$ associated 
with the real $1-$dimensional representation of Spin$(1,1)$ with spin $\pm \frac{1}{2}.$
Then we get (as in [page 587]\cite{Del}) a natural Dirac type operator $D_h$ acting on the sections of 
$(S_+\oplus S_-) \otimes \mathcal{E}_{h}$ and the map 
$$
s \in C^\infty( \frac{\C}{\Gamma^*} ;  \, (S_+\oplus S_-) \otimes \mathcal{E}_{h}) \rightarrow \int_{\frac{\C}{\Gamma^*}}  ( \overline{s} ; D_h(s) )(y) d y =L_h(s)
$$ defines a {\it free} Lagrangian $L_h$  on 
$\frac{\C}{\Gamma^*}$. 

\end{proof}

\begin{open question}  
 In  [page 277]\cite{Polchinski} (see also 
[page 558]\cite{Del}) it is proved  that the space of perturbative  renormalizable theories on $\R^4$ is an attractor 
for the renormalization group flow acting on an  infinite dimensional 
space  of lagrangians. Is it possible to develop  a (non pertubative) quantum field theory on the space 
 $ \frac{\C}{  \Gamma^*}$ in  which there should exist
 a natural infinite dimensional 
space  $\mathcal{I}$ of lagrangians such that 
  $\frac { \C \times {\rm T} \Gamma}{\Gamma}$ should be 
 an attractor  of $\mathcal{I}$ for the renormalization group flow?
\end{open question}

\subsection{ Further remarks and motivation of Section \ref{Section:curve} and 
\cite{El1}}
\label{subsection:motivation}$\;$

\medskip
Denote by $\mathcal{A}^j_{\mathcal{F}}(S(E_0) )$ 
($0 \leq j \leq 2$)  the set of  sections of the real vector bundle
$\wedge^j {\rm T}^* \mathcal{F}\rightarrow S(E_0) $ which are smooth along the leaves and continuous 
on $S(E_0)$. The metric $g$ of \eqref{disym} induces a metric $h_j(\cdot ,\cdot )$ on the bundle $\wedge^j {\rm T}^* \mathcal{F}$ ($0 \leq j \leq 2$),
we then denote by $L^2(S(E_0);  \wedge^j {\rm T}^* \mathcal{F})$ the   Hilbert completion of 
$\mathcal{A}^j_{\mathcal{F}}(S(E_0) ) $ with respect to the real scalar product:
\begin{equation}\label{scp}
\forall \omega, \omega^\prime \in \mathcal{A}^j_{\mathcal{F}}(S(E_0) ), \;
<\omega ; \omega^\prime> = \int_{S(E_0)} h_j(\omega , \omega^\prime )(\theta) d\mu(\theta).
\end{equation}
Let $d_{\mathcal{F}}^\dag$ denote the formal adjoint of the leafwise exterior derivative $d_{\mathcal{F}}$. Then one defines 
 the reduced $L^2-$leafwise real cohomology groups 
$H^j_{L^2, \tau}(S(E_0), \R)$, ($0\leq j\leq 2$) by
$$
H^j_{L^2, \tau}(S(E_0), \R)= \ker(\,  (d_{\mathcal{F}}^\dag )^*\cap (d_{\mathcal{F}}^\dag )^{* *}),\; 0 \leq j \leq 2
$$ where the unbounded operators $(d_{\mathcal{F}}^\dag )^*,  (d_{\mathcal{F}}^\dag )^{* *}$ both act 
on $L^2(S(E_0);  \wedge^j {\rm T}^* \mathcal{F})$. 

The following important result comes from \cite{D 2}
\begin{theorem} \label{orbit} {\item 1)} 
There is a natural bijection between the set of valuations 
$w$ of the function field $K(E_0)$ of $E_0$ and the set of primitive compact $\R-$orbits of 
$\phi^t$ on $S(E_0)$. 
It has the following property. If $w$ corresponds to $\gamma=\gamma_w$, then
$$ l(\gamma_w)= \log N(w).
$$
{\item 2)} Denote by $( \phi^t)^*_j $ the operator 
$( \phi^t)^* $ acting on $ H^j_{L^2, \tau}(S(E_0), \R)$ for $0\leq j \leq 2$. Then for 
any $\alpha \in C_c^\infty (\R ; \R)$ the  following equality holds: 
\begin{equation} \label{TE}
\sum_{j=0}^{j=2} (-1)^j \int_\R \, {\rm TR} ( \phi^t)^*_j \alpha(t)\, d t\,=\,
\sum_{\gamma_w} l(\gamma_w) \sum_{k\geq 1} \alpha({k l(\gamma_w)}) +
\sum_{\gamma_w} l(\gamma_w) \sum_{k\leq -1}  e^{kl(\gamma_w)} \alpha({k l(\gamma_w)}) 
\end{equation}
 where the right hand side coincides with the one of the explicit formula \eqref{EF2}. 

\end{theorem}
About part 2]. Deninger has identified the left handside of \eqref{EF2} with 
the (spectral) left handside of \eqref{TE} and has invoked  \eqref{EF2} to get
\eqref{TE}.

Now we make remarks about the structure of $(S(E_0), \mathcal{F})$ 
which will motivate the constructions and definitions of 
the Sections \ref{Section:curve} and \cite{El1}.

\smallskip
First we introduce carefully a natural transverse measure 
on $(S(E_0), \mathcal{F})$  and point out its important role.

Set $\mathcal{L}_{E_0}= \frac{ \C \times {\rm T} \Gamma}{\Gamma}$, this is a 
compact laminated space which is  foliated by its path-connected components. Any 
element $q^\nu \in q^\Z$ acts on $[z, \widehat{v}] \in \mathcal{L}_{E_0}$ 
by $q^\nu \cdot [z, \widehat{v}] = [\xi^{\nu} z, \xi^{\nu} \widehat{v}].$
The Haar measure $\mu_\xi$ of ${\rm T} \Gamma$ induces a transverse 
measure, still denoted $\mu_\xi$, of $\mathcal{L}_{E_0}.$ 
For any Borel transversal $T$ of $\mathcal{L}_{E_0}$ one has 
$\mu_\xi (q\cdot T) = q^{-1} \mu_\xi (T).$ 

Moreover, the metric $\widetilde{g}= (d x_1)^2 +  (d x_2)^2$ (where $z=x_1 +i x_2$) defines a 
leafwise metric on $\mathcal{L}_{E_0},$ let $\lambda_{\widetilde{g}} $ be the associated 
leafwise volume form. Then
 $\lambda_{\widetilde{g}} \mu_\xi$ defines a $q^\Z-$invariant measure of $\mathcal{L}_{E_0}.$

The leafwise metric $g$ in \eqref{disym}  of  $(S(E_0) , \mathcal{F})$ is 
defined by $g= x^{-1} \widetilde{g}$ and its associated leafwise volume form is given by
 $\lambda_g = x^{-1} dx_1 \wedge dx_2.$ 
 
 Now it is clear that   $\frac {\mu} { \lambda}={\mu_\xi d x}$ 
defines a transverse measure with associated 
Ruelle-Sullivan current $C(\frac {\mu} { \lambda}).$ We can pair sections of $\mathcal{A}^2_{\mathcal{F}}(S(E_0))$ with $C(\frac {\mu} { \lambda}),$ for instance the measure $\mu$ may be recovered by the formula:
$$
\forall f \in C^0( S(E_0)),\; ( f \lambda ; C(\frac {\mu} { \lambda} )  \,)= \int_{S(E_0)} f \, d \, \mu.
$$ Recall that $\mu$ is $\phi^t-$invariant. Moreover we observe that the scalar product 
\eqref{scp} 
may be recovered by the formula:
$$
\forall \omega, \omega^\prime \in \mathcal{A}^j_{\mathcal{F}}(S(E_0)),\; 
< \omega ; \omega^\prime > = ( \omega \cup \ast \omega^\prime ; C(\frac {\mu} { \lambda} )  \,)
$$ where $\ast$ denotes the leafwise Hodge star operator associated to $g$. 

Now we come to formula 
\eqref{TE}. As explained by Deninger, the dissymmetry of the coefficients of 
$\alpha({k l(\gamma)}) $ for $k\geq -1$ and $k\geq 1$  is due to 
property \eqref{disym} (see the remark following Corollary 1 of \cite{El}). 
We are going to propose a dynamical explanation, \` a la 
Guillemin-Sternberg, of this dissymmetry. Consider a 
point $(z_0, \widehat{v}_0, 1) \in S(E_0)$, with 
$  \widehat{v}_0 \in {\rm T} \Gamma$,  such that 
$\phi^{-\log q }[z_0, \widehat{v}_0, 1]=[z_0, \widehat{v}_0, q]=[z_0, \widehat{v}_0, 1].$
 Recall that by definition 
$(\xi^{-1} z_0, \xi^{-1} \widehat{v}_0, q^{-1} q ) \sim (z_0, \widehat{v}_0, q) $. So 
$[\xi^{-1} z_0, \xi^{-1}  \widehat{v}_0, 1] = [z_0, \widehat{v}_0, 1]$ and 
 there exists $\gamma \in \xi^{-1} \Gamma$ such that 
\begin{equation} \label{eq:fix}
\xi^{-1} z_0 = z_0 + \gamma,\quad \xi^{-1}  \widehat{v}_0 =  \widehat{v}_0 - \gamma.
\end{equation}
The operator $(\phi^t)^*$ acting on $ \mathcal{A}^j_{\mathcal{F}}(S(E_0))$ admits a 
Schwartz kernel defined by the formula:
$$
\forall \omega \in  \mathcal{A}^j_{\mathcal{F}}(S(E_0)), \; 
(\phi^t)^* ( \omega)(y)= \int_{S(E_0)} (D \phi^t)^* \delta_{\phi^t(y)=y^\prime}\,  
\omega(y^\prime)\,  d \mu(y^\prime).
$$ Consider a point $y=[ z, \widehat{v}, x]$ belonging to a small neighborhood 
 of $\{ \phi^t [z_0, \widehat{v}_0, 1],\, -\log q \leq t \leq 0\}.$ Then, with the previous notations, 
 one has: 
\begin{equation} \label{eq:local}
\phi^t(y)=(\xi^{-1} z-\gamma, \xi^{-1}\widehat{v} + \gamma, q^{-1} x e^{-t}).
\end{equation}
 In the 
following Lemma we show basically  that the graph of the flow $\phi^t$ is transverse to 
the diagonal and compute $\delta_{\phi^t(y)=y}.$ 
\begin{lemma} \label{lem:trans} {\item 1)} $ z \in \C \rightarrow \xi^{-1} z-\gamma-z$ and 
$\widehat{v} \in V_\xi \Gamma \rightarrow \xi^{-1}\widehat{v} + \gamma -\widehat{v}$ 
are invertible and their jacobians are respectively given by:
$$
{\rm Jac} (\xi^{-1} z-\gamma-z )= | \xi^{-1}-1 |^2  ,\; 
{\rm Jac} ( \xi^{-1}\widehat{v} + \gamma -\widehat{v})= {q} .
$$
{\item 2)}  Let $V$ be an open neighborhood of $(z_0, \widehat{v}_0)$, set:
$$
U=\{ (z, \widehat{v}, e^{-s})/ \, s\in ]-\log q, 0],\, (s,\widehat{v})\in V\}.
$$ Consider $\epsilon>0$ and  $V$ small enough so that 
 $t\in [-\log q, 0] \rightarrow (z_0, \widehat{v}_0, e^{- t})$ is the only 
 closed orbit of $\phi^t$ contained  in  $U$ with length in $]-\epsilon -\log q, \epsilon -\log q[.$
Then one has the following equality as a distribution on 
$U\times ]-\epsilon -\log q, \epsilon -\log q[:$
$$
\delta_{\phi^t(y)=y} = \frac{1}{| \xi^{-1}-1 |^2}\delta_{z-z_0}\otimes  \frac{1}{q} \delta_{\widehat{v}-\widehat{v}_0}
\otimes \,\delta_{t+ \log q}.
$$
\end{lemma}
\begin{proof} {\item1) } We prove only the second equality. Recall 
that ${\rm T} \Gamma$ is an open compact 
subset of  $V_\xi \Gamma.$ Then, since 
$\widehat{v} \rightarrow \widehat{v} - \xi \widehat{v}$ defines 
an automorphism of ${\rm T} \Gamma$ whose 
inverse is $\widehat{v} \rightarrow \sum_{n \in \N} \xi^n \widehat{v},$ one has 
${\rm Jac}\, ( \widehat{v} - \xi \widehat{v})=1.$  Now recall that  the proof 
of Lemma \ref{Haar} shows that $\mu_\xi ( \xi {\rm T} \Gamma) = 
\frac{1}{q} \mu_\xi (  {\rm T} \Gamma)$ so that 
$ {\rm Jac}\, ( \xi \widehat{v}) =\frac{1}{q}.$ By combining the last two 
equalities for Jac, one gets: 

$ {\rm Jac} ( \xi^{-1}\widehat{v} + \gamma -\widehat{v})= {q}.$
{\item2) } Using the change of variable formula 
for $\int d\mu_\xi$ and the equality $\xi^{-1} \widehat{v}_0 +\gamma-\widehat{v}_0=0$,
one sees that for 
$ \widehat{v}$ close to $ \widehat{v}_0$ one has 
$$
 \delta_{ \xi^{-1}\widehat{v} + \gamma -\widehat{v}}= 
\frac{1}{{\rm Jac} ( \xi^{-1}\widehat{v} + \gamma -\widehat{v}) }
\delta_{\widehat{v} -\widehat{v}_0}.
$$ Then a computation using  \eqref{eq:fix} and \eqref{eq:local} shows (see also [Section IV] \cite{C 2}),
that  for $y=[z,\widehat{v},x] \in U$ and $t\in ]-\epsilon -\log q, \epsilon -\log q[$ one has:
$$
\delta_{\phi^t(y)=y} = \frac{1}{ {\rm Jac} (\xi^{-1} z-\gamma-z )}\delta_{z-z_0}\otimes  \frac{1}{ {\rm Jac} ( \xi^{-1}\widehat{v} + \gamma -\widehat{v}) } \delta_{\widehat{v}-\widehat{v}_0}
\otimes   \,\delta_{t+ \log q}.
$$ By combining 1) with this equality  one gets the result.
\end{proof}

Recall now that $d \mu(y) = d x_1 d x_2 \otimes \mu_\xi \otimes \frac{d x }{x}.$ 
The formula $|\int_0^{-\log q} \frac{ d \,e^{-s}}{e^{-s}}|= \log q$ and  Lemma \ref{lem:trans}. 2] show that for $t$ close to $-\log q$ 
the distributional trace 
$$\int_{S(E_0)}  {\rm Tr}\, (D \phi^t)^* \delta_{\phi^t(y)=y}\,  
 d \mu(y)
$$ is well defined  (near $-\log q$) and equal to:
$$
\log q \sum_{\gamma_w,\, l(\gamma_w )= \log q}  \frac{1}{q} 
\delta_{-l(\gamma_w)}
$$ where $\gamma_w$ runs over the set of closed orbits 
of $\phi^t$ of length $l(\gamma_w)= \log q.$
Since $ {\rm Jac} ( \xi \widehat{v} + \gamma -\widehat{v})=1$ a similar argument shows that 
for $t$ close to $\log q$ 
the distributional trace 
$$\int_{S(E_0)}  {\rm Tr}\, (D \phi^t)^* \delta_{\phi^t(y)=y}\,  
 d \mu(y)
$$ is well defined  (near $\log q$) and equal to:
$$
\log q \sum_{\gamma_w, \,l(\gamma_w )= \log q}  
\delta_{l(\gamma_w)}.
$$ These observations motivate the definition of a transversally $p-$adic 
  complex laminated space that we shall introduce in 
Section \ref{Section:curve}. Deninger pointed out to us  that the condition \eqref{disym}
$(\phi^t)^*g=e^t g$ was probably too strong for being generalized. That is why 
in Proposition \ref{Wilson}.2] we shall replace it by 
$(\phi^t)^*[\lambda_g]=e^t[\lambda_g ]$ where $[\lambda_g ]$ denotes the reduced leafwise cohomology 
class of the leafwise kaehler form associated to $g$ (see Serre \cite{Serre} and 
Deninger-Singhof \cite{D-S}).

\smallskip
 Now let $Y$ be a smooth projective absolutely irreducible curve over $\F_q$ 
admitting a rational point. As pointed out in \cite{El}, one would like to associate to 
$Y$ a Riemannian laminated foliated space 
$( \frac{ \mathcal{L}_Y\times \R^{+*}}{q^\Z}, \mathcal{F}, g, \phi^t)$ for which the analogous version of 
Theorem \ref{orbit} should hold. Recall now Hurwitz's formula for 
a non constant morphism $\psi: C_1\rightarrow C_1$ where $C_1$ is a smooth projective curve 
of genus $g_1$ over a field of characteristic zero. One has:
$$
2 g_1- 2= ({\rm deg}\, \psi ) (2 g_1- 2) + \sum_{P \in C_1} (e_\psi(P) -1)
$$
where the strictly positive integers $e_\psi(P)$ are all equal to one except for a finite number of them. If it could be possible to lift the Frobenius morphism of 
$Y$ in characteristic zero, one should get a morphism $\psi$ of degree $>1$ which is not 
possible.
Thus, unlike the case of an elliptic curve,  we do not consider for $Y$ a flow of the simple
form $\phi^t ([l,x])= [l, x e^{-t}]$ but a priori a (more general) renormalization group flow.
We shall carry out the detailed constructions  in Section \ref{Section:curve}.

\section{The zeta function of a
 compact Riemannian foliated transversally $p-$adic laminated space  $(S=\frac {\mathcal{L} \times \R^{+*}}{q^\Z}, \mathcal{F}, g, \phi^t)$ 
 \label{Section:curve}}  $\;$

We introduce the following  definition which is a particular case of a notion due to Sullivan \cite{Sullivan}:
\begin{definition}  \label{def:trans} A compact complex $1-$dimensional  transversally $p-$adic laminated space $\mathcal{L}$
is a compact topological space satisfying the following property.
There exists $m \in \N^*$ such that $\mathcal{L}$ admits a finite open cover 
$\cup_{i=1}^N A_i $ by local charts: 
$ f_i: A_i \rightarrow  U_i \times \Z_p^m $ where $U_i$  is an open disc
of $\C$  such that each transition map
$f_i\circ f_j^{-1}$ is of the form 
\begin{align} & 
f_j(A_i\cap A_j)  \rightarrow U_i \times \Z_p^m
 \\&
(z, \theta) \rightarrow (H(z), G(\theta) )\,=\, f_i\circ f_j^{-1}(z,\theta)
\end{align} where $f_j(A_i\cap A_j)= \Omega_{i,j} \times \Z_p^m$ and
{\item 1)}  $z\rightarrow H(z)$ is holomorphic on its domain of definition $\Omega_{i,j}$

{\item 2)}  There exists $M\in {\rm GL}_m(\Z_p)$ and $B \in \Z_p^m$ such that 
for any $\theta \in \Z_p^m$ one has
$G(\theta)= M \theta + B.$  

Then the path connected components  $\mathcal{C}$
of $\mathcal{L}$ carry a structure of $1-$dimensional complex manifold and constitute 
the leaves of a foliation. For any $l \in \mathcal{L}$ we denote by 
 ${\rm T}_l \mathcal{L}$ the tangent space at $l$ 
of the path connected component of $\mathcal{L}$ containing $l.$ 

We may and shall assume {\bf in the sequel} that the charts $(A_j,f_j)$ are defined 
on open subsets $A^\prime_j$ such that 
$\overline{A_j}\subset A^\prime_j$.
\end{definition}
The space $\frac {\C \times V_\xi \Gamma} {V}$ of Section 2 is an example 
of such an $\mathcal{L}$.
\begin{definition}  \label{def:f} Let $\mathcal{L}$ be a laminated space as in Definition 
\ref{def:trans} and $j\in \{0,1,2\}.$ Denote by $\mathcal{A}^j_{\mathcal{F}}(\mathcal{L})$ 
the set of sections $\omega$ of the complex vector bundle 
$ \wedge^j {\rm T}^* \mathcal{L} \otimes \C \rightarrow \mathcal{L}$ 
which are smooth along the leaves and continuous globally in the following sense. 
In any chart $(A_i,f_i)$ as above, $(f_i^{-1})^*\omega$ is a finite sum of terms 
$a(z,\theta) (d x_1)^\alpha \wedge (d x_2)^\beta$ where 
$z=x_1 + \sqrt{-1} x_2, \alpha + \beta = j,$
$$
\forall \theta \in \Z_p^m,\, z=x_1 + \sqrt{-1} x_2 \rightarrow a(z,\theta) 
$$ belongs to $C^\infty(U_i, \C)$ and for any $r_1, r_2 \in \N,$ 
$(z,\theta) \rightarrow \partial^{r_1}_{x_1} \partial^{r_2}_{x_2} a(z,\theta)$ is continuous 
on $U_i \times \Z_p^m.$
\end{definition}

\smallskip Now we consider a $1-$dimensional complex compact transversally $p-$adic 
laminated space 
$\mathcal{L}.$
In the sequel, we shall 
 assume that  $\mathcal{L}$ satisfies the following four properties which are stated with the notations of Definition \ref{def:trans} and  motivated by 
subsection \ref{subsection:motivation}:

(i) 
The group $q^\Z$ acts leafwise holomorphically on  $\mathcal{L}$ 
in the following sense, $q^m \cdot l$ denoting the action of $q^m \in q^\Z$ on $ l \in \mathcal{L}.$ Let $y\in \mathcal{L}$, assume that $y \in A_j$ 
and $q\cdot y \in A_i$ (cf Definition \ref{def:trans}). Then: 
$$
\forall (z,\theta)  \in f_j(q^{-1}(A_i)\cap A_j),\; f_i \circ q \circ f_j^{-1}(z, \theta) \,=\,
(Z(z) , L(\theta)=M_q \theta + b_q)
$$ on its domain of definition, where $Z(\cdot) $ is holomorphic,  $b_q \in \Z_p^m$, and

\noindent $M_q \in {\rm M}_m(\Z_p) \cap {\rm GL}_m(\Q_p)$ is such that
${\rm Jac} \, M_q=\frac{1}{q} .$ 




 (ii) There exists a smooth  leafwise kaehler metric $\widetilde{g}$ along the leaves of $\mathcal{L}$ with the following two properties. 
First, for each foliation  chart $f_i: A_i \rightarrow U_i \times \Z_p^m$ as 
in Definition \ref{def:trans}, $(f_i^{-1})^*(\widetilde{g})$ does not depend 
on $\theta.$ Second,  for any $l \in \mathcal{L}$ 
one has:
$$
\forall u \in {\rm T}_l \, \mathcal{L},\; 
\widetilde{g}(q_*(u))=q\, \widetilde{g} (u).
$$ We denote by $\lambda_{\widetilde{g}}$ the leafwise Kaehler form associated with $\widetilde{g}$: $ \widetilde{g}(u,v)= \lambda_{\widetilde{g}}( u, J v).$

Since the leaves of $\mathcal{L}$  are oriented by their complex structure 
there is a leafwise Hodge star associated to $\widetilde{g}$.
 
(iii) There exists a smooth family $(\psi^t_x)_{t\in \R, x>0}$ of diffeomorphisms of 
$\mathcal{L}$ acting leafwise holomorphically in the following sense.
Let $(t_0,x_0, y)\in \R\times \R^{+*} \times \mathcal{L}$, assume that $y \in A_j$ 
and $\psi^{t_0}_{x_0}( y )\in A_i.$ 
Then there exists a small neighborhood $I_0\times J_0 \times B_y \subset 
\R \times \R^{+*}\times  (U_j \times \Z_p^m) $ of $(t_0,x_0,f_j(y))$ such that:
$$
\forall (t, x, z,\theta)  \in I_0\times J_0\times B_y,\; f_i \circ \psi_{x}^{t} \circ f_j^{-1}(z, \theta) \,=\,
(R_x^{t}(z) , E(\theta))
$$ where the  $(R_x^{t}(\cdot))_{(t,x)\in I_0\times J_0} $ defines a smooth family
of holomorphic maps on their domain of definition and  there exists 
$M_1 \in {\rm GL}_m(\Z_p)   $ and $b_1 \in \Z_p^m$ such that 
$E(\theta)= M_1 \theta + b_1$ on its domain of definition.
Moreover, for  any 
$ (t,t^\prime, x) \in \R \times \R \times \R^{+*}, \, l \in \mathcal{L}$ one assumes that:
\begin{equation} \label{Reno}
\psi^t_{q x }(q\cdot l) = q \cdot \psi^t_x(l),\;   \psi^{t^\prime}_{x e^{-t} }\left( \psi^t_x(l)\right) = \psi_x^{t+t^\prime}(l), 
\end{equation} 
 Notice that the equation \eqref{Reno}  implies
that  $\psi^0_x = {\rm Id},$  that each $\psi^t_x$ preserves any  path connected component of $\mathcal{L}$ and 
induces a diffeomorphism of it.

Moreover, let $k\in \N^*$ and $l_0\in A_j$ be such that 
$(\psi^{\log q}_q \circ q)^k (l_0)=l_0.$ Then one can write:
$$
\forall (z,\theta) \in U_j^\prime\times \Z_p^m,\;\, 
f_j\circ (\psi^{\log q}_q \circ q)^k \circ f_j^{-1} (z,\theta)= 
(Z_k(z), Q_k \theta + b_k)
$$ where $Z_k$ is holomorphic on an open neighborhood 
$U_j^\prime \subset U_j$ of $f_j(l_0)$, $b_k\in \Z_p^m$, and 
$Q_k \in {\rm M}_m(\Z_p)$ satisfies:
\begin{equation} \label{eq:jac}
{\rm Jac}(Q_k)=q^{-k},\; {\rm Jac}({\rm Id}-Q_k)= 1.
\end{equation}

(iv) Let $\omega \in C^0(\mathcal{L}, \wedge^1 {\rm T}^* \mathcal{L})$ 
be a leafwise $C^\infty$ harmonic $1-$form. Then for each chart 
$(A_i,f_i)$ of Definition \ref{def:trans},  $(f_i^{-1})^*\omega(z,\theta)$ does not depend 
on $\theta.$ For an intrinsic equivalent statement see Theorem \ref{Thm:intrinsic}. 

\medskip

\noindent {\bf Remark} The assumption (iv) is satisfied in the case of the space 
$\frac{\C \times V_\xi \Gamma}{V}$ of Section 2. Louis Boutet de Monvel 
has showed us an example (of the form 
$\frac{\mathcal{H}\times \Z_p^{2 g}}{ \Gamma}$) with hyperbolic leaves where it is not satisfied.
Dennis Sullivan has told us  that, in some sense, most of the examples of $\mathcal{L}$  with hyperbolic leaves
of Definition \ref{def:trans} do not satisfy this assumption. 
Bertrand Deroin has showed \cite{Deroin} that for a space $\mathcal{L}$ 
with hyperbolic leaves, the space of leafwise holomorphic 
 quadratic differentials is infinite dimensional. But his proof does not show 
that Assumption (iv) is not satisfied.
Notice that there are examples of spaces $\mathcal{L}$ which have 
hyperbolic and parabolic leaves and do not seem to be defined by an 
inverse limit of Riemann surfaces (eg the one described by E. Ghys 
[Section 6.4]\cite{Ghys} following an idea of R. Kenyon).

\begin{lemma} \label{lem:transv}{\item 1]} The data of the Haar measure $\mu_{\Z_p^m}$ in each 
local chart $(A_i,f_i)$ of Definition \ref{def:trans} induce a transverse measure 
$\mu_{\mathcal{L}}$ on $\mathcal{L}.$ Then 
$\lambda_{\widetilde{g}} \mu_{\mathcal{L}}$ defines a measure on $\mathcal{L}.$ Moreover,
 for any Borel transversal $T$ of 
$(\mathcal{L} , \mathcal{F})$ one has $\mu_{\mathcal{L}} (q\cdot T)= 
\frac{1}{q}  \mu_{\mathcal{L}}.$
{\item2]}
$$
\forall (t,x) \in \R \times \R^{+*},\;
(\psi^t_x)^*(\mu_{\mathcal{L}})\,= \,\mu_{\mathcal{L}}.
$$

\end{lemma} 
\begin{proof}1]. The fact that $\mu_{\mathcal{L}}$ is well defined as a transverse 
measure is a consequence of Definition \ref{def:trans}.  The leaves are two dimensional 
so $\lambda_{\widetilde{g}}$ defines the leafwise Riemannian volume form 
associated to $\widetilde{g}$ on $(\mathcal{L}, \, \mathcal{F}).$
Then by definition, $\lambda_{\widetilde{g}} \mu_{\mathcal{L}}$ defines a measure on $\mathcal{L}.$
The equality $\mu_{\mathcal{L}} (q\cdot T)= 
\frac{1}{q}  \mu_{\mathcal{L}}$ is a consequence of Assumption (i) ( 
${\rm Jac}\, M_q= \frac{1}{q}$).

\noindent 2]. This is a consequence of the equality ${\rm Jac} (M_1)=1$ in 
Property (iii).
\end{proof} 

\bigskip 
Consider the foliated laminated space 
$(S= \frac{ \mathcal{L} \times \R^{+*}}{q^\Z}, \mathcal{F}) $ whose leaves are 
 the sets $\Pi( \mathcal{C} \times \{x\}) $  where $\Pi:  \mathcal{L} \times \R^{+*} \rightarrow  S $ denote the projection and $\mathcal{C} $ runs over the set 
of path connected components of $ \mathcal{L} $.
 Let $[l,x_0] \in S$, assume that $l\in A_i$ with the notations 
of Definition \ref{def:trans}. There exists $\epsilon >0$ such that 
$ A_i \times ]x_0 -\epsilon, x_0+\epsilon[$ may be identified with an open subset of $S.$ Then
$$f_i \times {\rm Id}_{]x_0 -\epsilon, x_0+\epsilon[}:\, A_i \times ]x_0 -\epsilon, x_0+\epsilon[\rightarrow U_i\times \Z_p^m \times]x_0 -\epsilon, x_0+\epsilon[
$$ defines a foliated chart of $S$ centered at $[l,x_0].$

\begin{definition} \label{def:form} {\item 1]}
Denote by $\mathcal{A}^j_\mathcal{F}(S)$ ($0\leq j \leq 2 $) 
the set of sections $\omega$ of the complex vector bundle 

\noindent $\wedge^j {\rm T}^* \mathcal{F} \otimes \C \rightarrow S$ which are smooth along the leaves and continuous globally in the following sense. In any 
foliation chart $f_i \times {\rm Id}_{]x_0 -\epsilon, x_0+\epsilon[}$ as above,
$((f_i \times {\rm Id}_{]x_0 -\epsilon, x_0+\epsilon[})^{-1})^*\omega$ is  a finite sum of terms $a(z,\theta,x) (d x_1)^\alpha \wedge (d x_2)^\beta$ where $z=x_1 + \sqrt{-1} x_2,\,  \alpha + \beta=j$ such that:
$$\forall \theta \in \Z_p^m,\,
(z=x_1 + \sqrt{-1} x_2,x)  \rightarrow a(z, \theta,x)\,\in C^\infty(U_i \times ]x_0 -\epsilon, x_0+\epsilon[, \R),
$$  and for any $r_1, r_2, r_3 \in \N$, $(z,\theta,x) \rightarrow \partial_{x_1}^{r_1} \partial_{x_2}^{r_2} 
\partial_x^{r_3} a(z, \theta,x)$ is continuous. 
{\item 2]} Denote, for $j\in \{0,1,2\}$, by 
$\overline{H}^j_\mathcal{F}(S)$ the reduced leafwise cohomology 
$\frac{\ker d_{\mathcal{F}}}{ \overline{{\rm Im } d_{\mathcal{F}}}}.$
\end{definition}

\bigskip One then gets the following:
\begin{proposition}  \label{Wilson} Denote by $\Pi$ the projection map 
$\Pi: \mathcal{L} \times \R^{+*}\rightarrow\frac{ \mathcal{L} \times \R^{+*}}{q^\Z} =S$ where 
$q^\Z$ acts diagonally.
Consider the foliated laminated space 
$(S= \frac{ \mathcal{L} \times \R^{+*}}{q^\Z}, \mathcal{F}) $ whose leaves are 
 the sets $\Pi( \mathcal{C} \times \{x\}) $  where $\mathcal{C} $ runs over the set 
of path connected components of $ \mathcal{L} $.
{ \item 1]} One defines a flow $\phi^t$ acting on $S$ by setting for any 
$(l,x) \in  \mathcal{L} \times \R^{+*}$
$$
\phi^t([l,x]) = ( [\psi^t_x(l) , x e^{-t}]),
$$ where $[l,x] $ denotes the class of $(l,x)$ in $S.$
{ \item 2]} The metrics $ x^{-1}\widetilde{g} $ 
on 
${\rm T}_{[l,x]} \mathcal{F}$, define a leafwise kaehler metric $g =( x^{-1}\widetilde{g})_{x\in \R^{+*}} $ on 
$(S,\mathcal{F}).$ Moreover, $\lambda_g= x^{-1}\lambda_{\widetilde{g}} $ defines a 
leafwise kaehler form on $(\mathcal{S} , \mathcal{F})$ ($\lambda_{\widetilde{g}} $ is defined 
in Assumption (ii) of this Section).
For any $t\in \R,$ $(\phi^t)^* [ \lambda_g] = e^t\, [ \lambda_g]$ where $[ \lambda_g]$ 
denotes the induced class in $\overline{H}^2_\mathcal{F}(S).$ 

{ \item 3]}  $\lambda_{{g}}$ is also the leawise Riemannian volume form 
associated to ${g}$ on $(S,\mathcal{F})$. Moreover, 
$ \mu_{\mathcal{L}}\,  {d x } $ defines a transverse measure denoted $\Lambda$ on 
$(S, \mathcal{F})$, and
$\mu= \lambda_{{g}} \,\mu_{\mathcal{L}}\, {d x } $ induces a  measure 
on $S.$

{ \item 4]} Let $W( \mathcal{L}, \mathcal{F})$ (resp. $W( S, \mathcal{F})$) denotes the von Neumann 
algebra of the foliation $( \mathcal{L}, \mathcal{F}) $ (resp. $( S, \mathcal{F})$) associated 
with the measure $\lambda_{{\widetilde{g}}}  \,\mu_{\mathcal{L}} $ (resp. $\lambda_{{g}} \,\mu_{\mathcal{L}} \,d x   $).
 Assume that    $W( \mathcal{L}, \mathcal{F})$ is a factor. Then the flow $\phi^t$ induces an action, denoted $ ( \phi^t)^*$, on
 $W( S, \mathcal{F})$ by 
$$A\rightarrow ( \phi^t)^* \circ A \circ (\phi^{-t})^*
$$ where 
$A=(A_l)_{l\in S/\mathcal{F}}\in W( S, \mathcal{F})$ is a random operator.
The cross-product von Neumann algebra $W( S, \mathcal{F})\rtimes_{(\phi^t)^*} \R $ is a type ${\rm III}_{\frac{1}{q}}-$factor.
\end{proposition}
\noindent {\bf Remark} The algebra $W( S, \mathcal{F})\rtimes_{(\phi^t)^*} \R $ represents morally
the set of measurable functions on  the noncommutative space 
of the orbits of $\phi^t$ (ie the closed points). As explained in [Section 3]\cite{El} this 
matches with Connes's approach to the zeta function of a function field.
\begin{proof} 1] This assertion follows from \eqref{Reno}.

\noindent 2] The fact that $ x^{-1} \widetilde{g}$ defines a leafwise metric on 
${\rm T}_{[l,x]} \mathcal{F}$ is a consequence of Assumption (ii) above. 
Since $\psi^0_x= {\rm Id}$ and $q^*\lambda_{\widetilde{g}}=q \lambda_{\widetilde{g}},$ one checks that one can write
\begin{equation} \label{Eq:kaehler}
(\psi^t_x)^* ( \lambda_{\widetilde{g}})=  \lambda_{\widetilde{g}}+ x e^{-t} \lim_{k\rightarrow + \infty} 
d_\mathcal{F} B^{t,x}
\end{equation} 
where $(B^{t,x})_{t\in \R, x>0}$ is a smooth family of elements of 
$\mathcal{A}^1_{\mathcal{F}}(\mathcal{L})$ such that $ B^{t ,  x}= q^*  B^{t, qx}.$ 
See the proof of Lemma \ref{lem:cte} for a similar argument.
The fact 
that $(\phi^t)^* [ \lambda_{g}]  = e^t\,  [\lambda_{g} ]$ is then a direct consequence of \eqref{Eq:kaehler}.

\noindent 3] The fact that 
  $$
\mu=\lambda_{g} \,\mu_{\mathcal{L}}\,  {d x } =\lambda_{\widetilde{g}} \,\mu_{\mathcal{L}}\,  \frac {d x } {x}
$$ induces a   measure 
on $(S, \mathcal{F})$ is a consequence of Lemma \ref{lem:transv}.1]   and  of Assumption (ii) 
($q^*\widetilde{g}= q \widetilde{g}$). It is clear that the fraction
$$
\frac{\lambda_g \,\mu_{\mathcal{L}}\,  d x }{\lambda_{g}}=
 \mu_{\mathcal{L}}\,  {d x } 
$$ defines a transverse measure on 
$(S, \mathcal{F}).$

\noindent 4] Let us first check 
 that $W( \mathcal{L}, \mathcal{F})$ is a type ${\rm II}_\infty-$factor. 
Denote by $\tau$ the trace on $W( \mathcal{L}, \mathcal{F})$ induced by
 the transverse measure $\mu_{ \mathcal{L}}.$ 
Consider a 
chart $(A_i,f_i)$ as in Definition \ref{def:trans}. Consider 
a finite orthonormal family $h_j \in C^\infty_0(U_i, \R)$
($1\leq j\leq j_0$). So:
$$ \forall j,l \in \{0, \ldots, j_0\}\; 
 {\rm with}\,  j\not=l,\; \int_{U_i} h_j^2 \lambda_{\widetilde{g}}=1,\; 
\int_{U_i} h_j h_l \lambda_{\widetilde{g}}=0.
$$ For $k\in \N$, denote by $\pi(k,j_0) \in W( \mathcal{L}, \mathcal{F}) $ the projection 
defined by $$
f_i^*\left(  (\sum_{j=1}^{j_0} h_j \langle \cdot,  h_j \rangle ) 1_{p^k \Z_p^m} \right).
$$ It is clear, as $k$ and $j_0$ vary,  that the $\tau(\pi(k,j_0) )$ do not belong to a ladder of 
$\R^+$ and can be arbitrarily large. 
Therefore  $W( \mathcal{L}, \mathcal{F})$ is 
a type ${\rm II}_\infty-$factor. This implies easily that 
$W({S}, \mathcal{F})$ is a type  ${\rm II}_\infty-$von Neumann algebra 
whose center is $L^\infty( \frac{R^{+*}}{q^\Z}).$
Now 
denote by $\tau_\Lambda$ the trace on $W({S}, \mathcal{F})$ induced by 
$\Lambda$ 
(See [Section 3]{\cite{El}} for details). Then, Lemma \ref{lem:transv} and the definition of $\phi^t$ allow to see 
that $\tau_\Lambda \circ (\phi^t)^* = e^{-t} \tau_\Lambda.$ 
 Using Connes's results [pages 494 and 495]\cite{C 1}, one then checks that $W({S}, \mathcal{F})\rtimes_{(\phi^t)^*} \R$ is of type {\rm III}. Moreover, one has:
$$
S\left( W({S}, \mathcal{F})\rtimes_{(\phi^t)^*} \R\right) \cap \R^{+*}\,=\, \{ \lambda>0/\, (\phi^{\log \lambda}))^* = {\rm Id}\,\}
= q^\Z.
$$ Then [page 473]\cite{C 1} implies that 
 $W({S}, \mathcal{F})\rtimes_{(\phi^t)^*} \R$ is of type 
${\rm III}_{\frac{1}{q}}.$ 
\end{proof}
\medskip
We think of $\mathcal{L}$ as a set of renormalizable quantum field theories, if we write $\psi^t_x= R_{x , x \,e^{-}t}$ 
then we recognize the general scheme of the method of the  renormalization group flow \`a la Wilson ([page 554]\cite{Del}). The condition 
$\psi^t_{q x }(q\cdot l) = q \cdot \psi^t_x(l) $ in \eqref{Reno} means that 
$q$ induces an action on $\mathcal{L}$ which commutes with the renormalization group flow 
up to rescaling: $R_{q x , q x \,e^{-t}}(q \cdot l) = q\cdot R_{x , x \,e^{-t}}(l)$ for 
any $l \in \mathcal{L}$. 



\bigskip Let $j\in \{0,1,2\}.$ 
Denote by $\mathcal{H}_\tau^j \subset \mathcal{A}^j_\mathcal{F}(S)$  
(see Definition \ref{def:form}) the subspace of 
complex leafwise harmonic forms and by $\pi^j_\tau$ the orthogonal  projection 
onto $\mathcal{H}_\tau^j \subset \mathcal{A}^j_\mathcal{F}(S)$, see Theorem \ref{Theorem:Hodge} and Proposition \ref{Prop:cte} for more on $\pi^j_\tau$. 
Now, denote by $C(\Lambda) $ the Ruelle-Sullivan current associated to the 
transverse measure $\Lambda = \mu_{\mathcal{L}}  d x $ of 
$(S, \mathcal{F})$. 
Then,
one defines a scalar product 
on $ \mathcal{H}_\tau^j $ by the formula 
$< \omega ; \omega^\prime> = ( \omega \cup \ast \overline{ \omega^\prime} ; C(\Lambda) )$ 
where $\ast$ denotes the leafwise Hodge star of the metric $g.$
We 
denote by  $H_\tau^j$ 
the $L^2-$completion of $ \mathcal{H}_\tau^j $.
 The operator $\pi^j_\tau \circ (\phi^t)^*$ defines a one parameter group acting on 
$H^j_\tau$  (see Proposition \ref{Prop:cte} for details). We
write $e^{t \Theta_j} = \pi^j_\tau \circ (\phi^t)^*$ where the infinitesimal generator $ \Theta_j$ 
defines an unbounded operator on the Hilbert space $H^j_\tau$.

 \begin{theorem} \label{Thm:trace}  Assume that the closed orbits $\gamma$ of the flow $\phi^t$ acting on $S= \frac{\mathcal{L}\times \R^{+*}}{q^\Z}$ are 
non degenerate.
Assume the four properties (i) to (iv),  that $W(\mathcal{L}, \mathcal{F})$ is a factor and 
that $\mathcal{L}$ has a dense leaf. Let $\alpha \in C^\infty_c(\R, \R)$.
{\item 1]} For each $j\in\{0,1,2\}$, $ \int_\R \alpha(t) e^{t \Theta_j}  d t$ 
acting on $H^j_\tau$ is trace-class and one has:
$$
 {\rm TR} \int_\R \alpha(t) e^{t \Theta_0}  d t \,=\,
\sum_{\nu \in \Z}   \int_\R \alpha(t) e^{\frac { 2 i \pi \nu t}{\log q} }  d t, \;
{\rm TR}  \int_\R \alpha(t) e^{t \Theta_2}  d t \,=\,
\sum_{\nu \in \Z}  \int_\R \alpha(t) e^{t+ \frac { 2 i \pi \nu t}{\log q} }  d t.
$$ Moreover, there exists a finite subset $\{ \rho_1,\ldots, \rho_{2 g}\} \subset \C$ such that
$$
{\rm TR} \int_\R \alpha(t) e^{t \Theta_1}  d t \,=\,
\sum_{j=1}^{2 g}\sum_{\nu \in \Z}   \int_\R \alpha(t) e^{t\, ( \rho_j+\frac {   2 i \pi \nu }{\log q}) }  d t
$$ 
{\item 2]} One has 
$$
\forall j \in \{1,\ldots, 2 g\},\quad \Re \rho_j= \frac{1}{2}.
$$
{\item 3]} One has
$$
\sum_{j=0}^{j=2} (-1)^j \, {\rm TR} \int_\R \alpha(t) e^{t \Theta_j}  d t\,=\,
$$
$$
\chi_\Lambda (\mathcal{F}) \alpha(0) + \sum_{\gamma} \sum_{k \geq 1} l(\gamma)
\left( \,e^{- k l(\gamma) } \alpha( -k l(\gamma) ) + \alpha( k l(\gamma) )\, \right)
$$ where $\chi_\Lambda (\mathcal{F})$ denotes Connes's $\Lambda-$Euler 
characteristic of $(S, \mathcal{F})$ (\cite{C 1}) and
$\gamma$ runs over the set of primitive closed orbits of $\phi^t$ 
and $l(\gamma)$ denotes the length of $\gamma.$
\end{theorem}


\noindent 
{\bf Remark.}  Define the Ruelle zeta function 
$$ \zeta_S(s)= \Pi_{\gamma} \,\frac{1}{1 -\exp\left( {-s \, l(\gamma)}\right)}
$$ where $\gamma$ runs over the set of primitive compact orbits of $\phi^t$ and 
$l(\gamma)$ denote its length. Deninger told us that 
Illies's result \cite{Illies}  and Theorem \ref{Thm:trace}  should imply that 
$\zeta_S(s)$ is an alternate product of regularized determinants. 
Then Theorem \ref{Thm:trace} should imply (along the lines of Deninger's formalism):
meromorphic extension of $\zeta_S(s)$ to 
$\C$,  functional equation $ s \leftrightarrow 1-s,$ Riemann hypothesis (or Weil's Theorem type) for $ \zeta_S(s).$

\begin{open question}
{\item 1]} Let $Y$ be a smooth projective absolutely irreducible curve over $\F_q$ admitting a rational point. Does there exist 
 a laminated foliated space $(S_Y=\frac {\mathcal{L_Y} \times \R^{+*}}{q^\Z}, \mathcal{F}, g, \phi^t)$ satisfying all the  assumptions of Proposition \ref{Wilson} and Theorem \ref{Thm:trace} 
and the following  
assumption:

(A) One has a natural bijection 
$w \mapsto \gamma_w$ between the set of closed points of $Y$ and 
the set of primitive closed orbits of $(S_Y, \phi^t)$ satisfying $\log N w= l(\gamma_w)$.

If the answer is yes, one should obtain, via Theorem \ref{Thm:trace}, a new proof 
of Weil's Theorem.
{\item 2]} Is it possible to interpret $\mathcal{L}_Y$ as an attractor 
of the renormalization group (semi-)flow acting on a suitable  infinite dimensional 
space of lagrangians? (compare with the end of \cite{D 2b} and [pages 30 and 559]\cite{Del}).
\end{open question}

Of course, when $Y$ is an elliptic curve over $\F_q,$ Deninger has 
shown (\cite{D 2}) that   the answer to part 1] of the previous open question
is yes. 

\section{Analytic results on  $\mathcal{L}$} $\;$


\subsection{Sobolev spaces on $\Z_p^m$} $\;$

\medskip
Recall that a character $\chi$ of $\Z_p^m=\oplus_{j=1}^m \Z_p$ may be written as:
$$
\chi= \oplus_{j=1}^m \sum_{l=0}^{n_j} \frac{a_{l,j}}{p^l}
$$ where the $a_{l,j}$ belong to $\{0,\ldots, p-1\}\subset \Z_p$ 
and $a_{n_j,j} \not=0$ if $n_j \geq 1.$ Notice 
that the $a_{l,j}$ are unique for $l\geq 1.$
If $\chi \not=0,$  we set 
$$
| \chi |= \max_{1\leq j\leq m,\, n_j\geq 1} | a_{n_j,j} |_p.
$$ If $\chi=0,$ we set $| 0 |=0.$
\begin{definition} \label{def:Delta_p} {\item 1]} Let $k\in \N$ and $u \in L^1(\Z_p^m, \, d \mu_{\Z_p^m}).$ We say 
that $u$ belongs to the Sobolev space 
$H^k(\Z_p^m)$ if the function
$$
\theta \rightarrow \sum_{\chi \in \widehat{\Z_p^m}} | \chi |^{k} \widehat{u}(\chi) \langle \chi, -\theta  
\rangle
$$ belongs to $L^2(\Z_p^m, \, d \mu_{\Z_p^m}).$ Here, $\widehat{u}(\chi)$ denotes the Fourier 
transform with the convention $\widehat{1_{\Z_p^m}}(0) =1.$
{\item 2]} The operator $\Delta_p$ defined by 
$$\Delta_p (u)(\theta)\,=\, 
 \sum_{\chi \in \widehat{\Z_p^m}} | \chi |^{2} \widehat{u}(\chi) \langle \chi, -\theta  
\rangle
$$ induces a bounded operator from $H^{k+2}(\Z_p^m)$ to $H^k(\Z_p^m)$ for any 
$k\in \N.$
\end{definition}
\subsection{Sobolev spaces and harmonic forms on $\mathcal{L}$} $\;$

\smallskip
The measure $\lambda_{\widetilde{g}} \mu_{\mathcal{L}}$ (see Lemma 
\ref{lem:transv}) and the leafwise metric $\widetilde{g} $   allow to define a scalar product 
on $C^0(\mathcal{L}, \wedge^j {\rm T}^* \mathcal{L})$ ($0\leq j \leq 2$). 
The axioms of Definition \ref{def:trans} (eg the fact that $M\in {\rm GL}_m(\Z_p)$) 
show that the following definition makes sense.
\begin{definition} \label{def:Delta} 
Let $(k,l) \in \N^2.$ We say that a function $u: \mathcal{L} \rightarrow \C$
belongs to the Sobolev space $H^{k,l}(\mathcal{L}, \C)$ if for any chart 
$(A_j,f_j)$ as in Definition \ref{def:trans} and any $\alpha\in \N^2$ with $|\alpha|\leq l$ the function:
$$
(z,\theta) \rightarrow  \sum_{\chi \in \widehat{\Z_p^m}} (1+ | \chi |^{k} )\,
\partial_{x_1}^{\alpha_1} \partial_{x_2}^{\alpha_2}\widehat{u\circ f_j^{-1}}^2(z,\chi) \langle \chi, -\theta  
\rangle
$$ belongs to $L^2(U_j\times \Z_p^m,\; d x_1 d x_2 \otimes d \mu_{\Z_p^m}).$ 
Here, $\widehat{\, }\,^2$ denotes the Fourier tranform with respect to the 
second variable (ie $\theta$).
Similarly one defines Sobolev spaces $H^{k,l}(\mathcal{L}, \wedge^* {\rm T}^* 
\mathcal{L})$ of differential forms.
\end{definition}

\begin{lemma} \label{lem:Delta} One defines a transverse $p-$adic Laplacian $\Delta_{p,T}$ acting on
 $\mathcal{A}^*_\mathcal{F}( \mathcal{L})$ (see Definition \ref{def:f})   by setting in each chart $(A_j,f_j)$ as in Definition \ref{def:trans}: 
$$
\forall u \in \mathcal{A}^*_\mathcal{F}( \mathcal{L}),\;
\Delta_{p,T}(u_{| A_j})\,= \,(f_j)^* \Delta_p (  (f_j^{-1})^* u_{| A_j}).
$$
\end{lemma}
\begin{proof} We prove the result for leafwise differential forms of degree one.
Consider $u \in \mathcal{A}^1_\mathcal{F}( \mathcal{L})$ having compact support 
in $A_i \cap A_j$ (with the notations of Definition \ref{def:trans}). We are going to show that 
$$
\left( f_i \circ f_j^{-1}\right)^*\Delta_p ( ( f_i^{-1})^* u )\,=\,\Delta_p (  ( f_j^{-1})^* u ).
$$ This will prove that  the operator $\Delta_{p,T} $ is intrinsically defined on $\mathcal{L}.$
From the assumptions of Definition \ref{def:trans} we deduce that 
$f_i(A_i\cap A_j)$ is of the form $ \Omega_{j,i}\times \Z_p^m.$ Recall 
that
$$
 f_i \circ f_j^{-1}(z,\theta)=(H(z), G(\theta)=M \theta + B)
$$ with $M\in {\rm GL}_m(\Z_p)$.

For any tangent vector $v \in {\rm T}_z \C$, we then have
$\left( f_i \circ f_j^{-1}\right)^*\Delta_p ( ( f_i^{-1})^* u )(z,\theta)=$
$$
\sum_{\chi \in \widehat{\Z_p^m}} | \chi |^2 
\int_{\Z_p^m}  ( ( f_i^{-1})^* u)(H(z),\xi)(D_zH(v))\, <\chi, \xi-G(\theta)> d \mu_{\Z_p^m}(\xi).
$$
In this integral we make the change of variable $\xi=G(\theta^\prime).$ We then have 
$$
< \chi , G(\theta^\prime) -G(\theta)>= < \chi , G(\theta^\prime-\theta)>= <^{t}G(\chi) , \theta^\prime-\theta >.
$$ Since
$M\in {\rm GL}_m(\Z_p)$ we observe that $\chi \rightarrow ^{t}G(\chi)$ defines a bijection 
of  $\widehat{\Z_p^m}$ satisfying $| ^{t}G(\chi) |= | \chi |.$ One then gets immediately that 
$$
\left( f_i \circ f_j^{-1}\right)^*\Delta_p ( ( f_i^{-1})^* u )\,=\,\Delta_p (  ( f_j^{-1})^* u ).
$$
\end{proof}
\begin{theorem} \label{Thm:intrinsic} {\item 1]} Assumption (iv) of Section \ref{Section:curve} is equivalent to the following 
(intrinsic) condition. Any leafwise $C^\infty$ harmonic $1-$form 
$\omega \in C^0(\mathcal{L}, \wedge^1 {\rm T}^* \mathcal{L})$ satisfies 
$\Delta_{p,T} \omega=0$.
 {\item 2]} 
Assume Assumptions from (i) to (iv) of Section \ref{Section:curve}.  Then the vector space $\mathcal{H}^1_{\mathcal{L}}$
of real harmonic $1-$forms $\omega \in C^0(\mathcal{L},\, \wedge^1 {\rm T}^* \mathcal{L})$ is of finite dimension $2 g$ where $g\in \N.$
\end{theorem}
\begin{proof} 
1] Left to the reader

\noindent 2]
Denote by $\Delta_{\mathcal{L}}$ the (essentially self-adjoint) leafwise signature 
laplacian associated with the metric $\widetilde{g}$ as in Assumption (ii) of Section \ref{Section:curve}. Notice that 
$ \Delta_{\mathcal{L}} $ is constant in $\theta$ in each chart $(A_i, f_i).$
Since $ \mathcal{L}$ is compact, an ellipticity argument shows that the operator 
$ \Delta_{\mathcal{L}} + \Delta_{p,T}$ is Fredholm 
from $\cap_{k+l=2} H^{k,l}(\mathcal{L}, \wedge^* {\rm T}^* 
\mathcal{L}) $ to $ L^2(\mathcal{L}, \wedge^* {\rm T}^* 
\mathcal{L})$. Indeed, one constructs a parametrix by considering in each 
chart $(A_j,f_j)$ the operator:
$$
\omega (z,\theta) \rightarrow \sum_{\chi \in \widehat{\Z_p^m}} 
\left( 1 + | \chi |^2 + \Delta_{\mathcal{L}} \right)^{-1} \widehat{\omega}^2(z,\chi) 
\langle \chi , -\theta\rangle.
$$
Now, observe that Assumption (iv) implies that 
$\mathcal{H}^1_{\mathcal{L}}\subset \ker (\Delta_{\mathcal{L}} + \Delta_{p,T}).$  Then
one gets that $\mathcal{H}^1_{\mathcal{L}}$ is finite dimensional.
Moreover,  the Hodge star $\star$ induces a complex structure on 
$\mathcal{H}^1_{\mathcal{L}}$ so this dimension is an even integer $2 g$.
\end{proof}
\begin{proposition} \label{Prop:basis}
The scalar product defined at the beginning of this Subsection induces 
a hermitian scalar product $ < ; >$ on  $\mathcal{H}^1_{\mathcal{L}}\otimes_\R \C$. There exists an orthonormal basis 
$(\omega_1,\ldots, \omega_{2 g})$ of $\mathcal{H}^1_{\mathcal{L}}\otimes_\R \C$ of eigenvectors 
for the action of $q$ on $\mathcal{H}^1_{\mathcal{L}}\otimes_\R \C.$
More precisely,  for each $j \in \{1,\ldots, 2 g \}$, there exists 
$\rho_j \in \C$ such that $\Re \rho_j = \frac{1}{2}$ and 
$q^*(\omega_j)= q^{\rho_j } \omega_j.$ 

\end{proposition}
\begin{proof} Lemma \ref{lem:transv}.1] and  Assumption  (ii) of Section \ref{Section:curve} 
(eg $q^*(\widetilde{g})= q \widetilde{g}$ ) allow to see that 
$q^*$ induces an operator acting on $\mathcal{H}^1_{\mathcal{L}}\otimes_\R \C$ 
of the form $ \sqrt{q} \, U$ where $U$ is unitary. This proves the result.

\end{proof}
\section{ Proof of Theorem \ref{Thm:trace}}$\;$

\subsection{Leafwise Hodge Decomposition and Heat operator}$\;$

We begin with describing a finite system of local foliated charts 
of $(S=\frac{\mathcal{L} \times \R^{+*}}{q^\Z}, \mathcal{F})$ with the notations of Definition \ref{def:trans}.

Set $I_1= ] q^{\frac{1}{2}}, q[$ and $I_2=]q^{-\frac{1}{10}} , q^{\frac{2}{3}}[.$ 
Set $\widetilde{A}_{j,2}=\{ [(q^k\cdot l, q^k x)]_{k \in \Z},\; 
(l,x) \in A_j \times I_2 \}$ and 
$F_{j,2}( [(q^k\cdot l, q^k x)]_{k \in \Z})=(f_j(l),x) \in U_j\times \Z_p^m \times I_2.$ 

Set also 
$\widetilde{A}_{i,1}=\{ [(q^k\cdot l, q^k x)]_{k \in \Z},\; 
(l,x) \in A_i\times I_1 \}$ and 

\noindent 
$F_{i,1}( [(q^k\cdot l, q^k x)]_{k \in \Z})=(f_i(l),x) \in U_i\times \Z_p^m \times I_1.$

Then, the $(\widetilde{A}_{i,1}, F_{i,1}),\, (\widetilde{A}_{j,2}, F_{j,2})$ 
($1\leq i, j \leq N$) define a finite open cover of foliated charts of 
$(S, \mathcal{F}).$ The transition maps are given in the following way.
If $(z,\theta,x)\in f_j(A_j \cap A_i) \times ]q^{\frac{1}{2}}, q^{\frac{2}{3}}[,$
one has 
$$
F_{i,1} \circ F_{j,2}^{-1}( z,\theta,x)= (f_i\circ f_j^{-1}(z,\theta),x).
$$ 

\noindent
If $(z,\theta,x)\in f_j(A_j \cap q^{-1}(A_i)) \times ]q^{-\frac{1}{10}}, 1[,$
one has 
$$
F_{i,1} \circ F_{j,2}^{-1}( z,\theta,x)=(f_i\circ q\circ  f_j^{-1}(z,\theta), q x).
$$

We can now state the following:
\begin{definition} Let $k\in \N$ and $j\in \{0,1,2\}.$ We say that an $L^2-$leafwise differential form 
$\omega \in L^2(S, \wedge^j T^*\mathcal{F})$ belongs to the 
space $H_{0,k}(S,\wedge^j T^*\mathcal{F} )$ if in any chart of the type $(\widetilde{A}_{i,1}, F_{i,1}),\, (\widetilde{A}_{j,2}, F_{j,2})$ as above,
$P \omega(z,\theta,x) \in L^2 $ for any differential operator of degree $\leq k$
 $P$ in the variables $(\Re z, \Im z, x).$   We set:
$$
H_{0,+\infty}(S,\wedge^j T^*\mathcal{F} )= \cap_{k\in \N}H_{0,k}(S,\wedge^j T^*\mathcal{F} ).
$$
\end{definition}

\begin{theorem} \label{Theorem:Hodge} Let $j\in \{0,1,2\}.$ 
{\item 1]}
We have the following leafwise Hodge decomposition:
$$
H_{0,+\infty}(S,\wedge^j T^*\mathcal{F} )= \mathcal{H}_{\tau,0,+\infty}^j 
\oplus^\perp \overline{ {\rm Im } \,d _{\mathcal{F}}} \oplus^\perp \overline {{\rm Im } \,\delta}
$$ where  $ \mathcal{H}_{\tau,0,+\infty}^j $ denotes the set of 
leafwise harmonic forms,  $\delta$ denotes the adjoint of the leafwise 
exterior derivative
$ d _{\mathcal{F}}$ of $S.$ The orthogonality 
is with respect to the $L^2-$scalar product. Moreover, one has 
$$
\mathcal{H}_{\tau,0,+\infty}^1= \mathcal{H}_{\tau,0,+\infty}^{1,0}  \oplus 
\mathcal{H}_{\tau,0,+\infty}^{0,1}
$$ where $\mathcal{H}_{\tau,0,+\infty}^{1,0}$ and $ \mathcal{H}_{\tau,0,+\infty}^{0,1}$
denote respectively the harmonic forms of type $(1,0)$ and $(0,1).$
{\item 2]}
 Denote by $\pi^j_\tau$ the orthogonal projection from $H_{0,+\infty}(S,\wedge^jT^*\mathcal{F} )$
onto $ \mathcal{H}_{\tau,0,+\infty}^j .$ Then $\pi^j_\tau \circ (\phi^t)^*$ defines 
a one parameter group acting on $ \mathcal{H}_{\tau,0,+\infty}^j .$
\end{theorem}
\begin{proof} {\item 1] } The foliated space $(S,\mathcal{F})$ is morally closed to a 
Riemannian foliation. One has just to adapt to our context the proof of Theorem 
1.1 of \cite{A-K1}. We leave the details to the reader.
{\item 2] } Denote, for $j\in \{0,1,2\},$  by $\overline{H}^j_{0,+\infty}$ the reduced leafwise cohomology group 
associated with $H_{0,+\infty}(S,\wedge^j T^*\mathcal{F} ).$
 Part 1] allows to identify $\mathcal{H}_{\tau,0,+\infty}^j$ with $\overline{H}^j_{0,+\infty}$
 on which $(\phi^t)^*$ acts naturally. One then gets easily the 
result.
\end{proof}

Now we introduce a particular partition of unity of $S$. Consider $\beta_1   \in C^{+\infty}_c(]q^{\frac{1}{2}}, q[ , \R)$  and 
$\beta_2 \in C^{+\infty}_c(] q^{\frac{-1}{10}}, q^{\frac{2}{3}}[  , \R)$  such 
that:
$$
\forall x \in [q^{\frac{1}{2}}, q^{\frac{2}{3}}],\; \beta_1(x)+ \beta_2(x) =1,\; 
\forall x \in [ q^{\frac{9}{10}}, q],\; \beta_1( x )+ \beta_2(q^{-1} x)=1, 
$$ 
$$
\forall x \in [1, q^{\frac{1}{2}}],\; \beta_2(x)=1,\; \forall x \in [q^{\frac{2}{3}}, q^{\frac{9}{10}}],\; \beta_1(x)=1.
$$

Next we set $\widetilde{V}_j= \{ [(q^k\cdot l, q^k x)]_{k \in \Z},\; l\in \mathcal{L}, x \in I_j\}$ for $j=1,2.$ If 
$ [(q^k\cdot l, q^k x)]_{k \in \Z} \in \widetilde{V}_j$ with 
$x\in I_j$, we set 
$$
\widetilde{F}_j( [(q^k\cdot l, q^k x)]_{k \in \Z}) = (l,x),\; 
\widetilde{\beta}_j( [(q^k\cdot l, q^k x)]_{k \in \Z})= \beta_j(x).
$$ By construction one has 
$ \widetilde{\beta}_1(y) + \widetilde{\beta}_2(y) \equiv 1$ on $S.$  Now we may state the:
\begin{definition}  \label{Def:contraction} If $\omega \in H_{0,k}(S,\wedge^j T^*\mathcal{F} )$ and 
$t\in \R^{+*},$ we set:
$$  
\mathcal{R}_t (\omega) = 
\sum_{j=1}^2 \left( \widetilde{F}_j^{-1} e^{-t \Delta_{p,T} }\widetilde{F}_j \right) ( \widetilde{\beta}_j \omega)
$$ where $  \Delta_{p,T}$ is defined in Lemma \ref{lem:Delta}.
\end{definition}
\noindent {\bf Remark.} Since the operator $\Delta_p$ of Definition \ref{def:Delta_p}   does not 
commute nicely with the multiplication by $p$ one cannot 
define a $p-$adic transversal Laplacian on $(S,\mathcal{F}).$ But the operator $\Delta_{p,T}$ 
allows to define  such a 
$p-$adic transversal Laplacian in each open subset $\widetilde{V}_1, \widetilde{V}_2$ 
of $S.$

\begin{proposition} \label{Prop:cte} 
{\item 1] } Let $\omega \in 
\mathcal{H}^j_{\tau,0,+\infty}.$ Then for any real $t>0,$ 
$\mathcal{R}_t (\omega) = \omega$ and,  $\omega$ is continuous on 
$S$ and  constant in $\theta$ in the above local charts $(\widetilde{A}_{i,1}, F_{i,1}),\, (\widetilde{A}_{j,2}, F_{j,2})$. Thus
$\mathcal{H}^j_{\tau,0,+\infty}= \mathcal{H}^j_{\tau}$ 
(see notation before Theorem \ref{Thm:trace}) and 
$$
\mathcal{H}^{1,0}_{\tau,0,+\infty}= \mathcal{H}^{1,0}_{\tau},\quad 
\mathcal{H}^{0,1}_{\tau,0,+\infty}= \mathcal{H}^{0,1}_{\tau}.
$$
{\item 2]} The operators $\pi^j_\tau \circ (\phi^t)^*$ define 
a one parameter group acting on  the Hilbert space ${H}^j_{\tau}$ of Theorem  \ref{Thm:trace}. 
\end{proposition}
\begin{proof} {\item 1]}  Assumption (ii) of Section \ref{Section:curve} states that the metric 
$\widetilde{g}$ is constant in $\theta$ in the chart $(A_i,f_i).$ Then, for each $t>0$, one checks
 that 
$ \mathcal{R}_t (\omega)$ is a 
smooth leafwise harmonic form which is 
 continuous on $S.$ Assumption (iv) of Section \ref{Section:curve} allows to see that  $ \mathcal{R}_t (\omega)$  is 
 constant in $\theta$ in the charts $(\widetilde{A}_{i,1}, F_{i,1}),\, (\widetilde{A}_{j,2}, F_{j,2})$.
 Letting $t$ go to $0^+$ one gets 
that $\omega$ is  constant in $\theta$ in these local charts and is equal to $ \mathcal{R}_t (\omega)$ for 
any $t>0.$
{\item 2]}  One checks that $(\phi^t)^*$ acts continuously on $ H_{0,+\infty}(S,\wedge^j T^*\mathcal{F} )$
with respect to the $L^2-$scalar product. The result then follows from
part 1] and  Theorem \ref{Theorem:Hodge}. 2].
\end{proof}

\subsection{Proof of Theorem \ref{Thm:trace}. 1]}$\;$

Recall that 
$e^{t \Theta_j}= \pi_\tau^j\circ (\phi^t)^*$. Actually the results of this Subsection 
 allow to (re)prove, using explicit computations,  that $\pi_\tau^j\circ (\phi^t)^*$ defines a one parameter 
group acting on $H^j_\tau.$ 

Recall that $\mathcal{L}$ has a dense leaf.
Then a continuous function on  $\mathcal{L}$ which is constant along the leaves is constant. Then one checks easily 
that $(x^{\frac{2 i \pi \nu}{\log q}})_{\nu \in \Z}$ is an orthonormal basis of 
$H^0_\tau \otimes_\R \C$ and that:
$$
 {\rm TR} \int_\R \alpha(t) e^{t \Theta_0}  d t \,=\,
\sum_{\nu \in \Z}   \int_\R \alpha(t) e^{\frac { 2 i \pi \nu t}{\log q} }  d t.
$$
In the same way, one checks that $(x^{\frac{2 i \pi \nu}{\log q}} 
x^{-1}\lambda_{\widetilde{g}})_{\nu \in \Z}$ is an orthonormal basis of 
$H^2_\tau \otimes_\R \C$ (recall that $\lambda_{\widetilde{g}}$ is defined in Assumption 
(ii) of Section \ref{Section:curve}).   Then one easily checks that:
$$
{\rm TR}  \int_\R \alpha(t) e^{t \Theta_2}  d t \,=\,
\sum_{\nu \in \Z}  \int_\R \alpha(t) e^{t+ \frac { 2 i \pi \nu t}{\log q} }  d t.
$$  
Next,   using  Proposition \ref{Prop:basis} and its notations, one checks that 
$$
(x^{\frac{2 i \pi \nu}{\log q}-\rho_j} 
\omega_j)_{j\in \{1,\ldots, 2 g\},\, \nu \in \Z}
$$ is an orthonormal basis of 
$H^1_\tau \otimes_\R \C.$ 
\begin{lemma} \label{lem:cte} Let $j\in \{1,\ldots, 2 g\}$ and $\nu \in \Z.$ Then 
one has 
$$(\phi^t)^* (x^{-\rho_j -\frac{2 i \pi \nu}{\log q} } \omega_j)= 
e^{t( \rho_j +\frac {2 i \pi \nu}{\log q}) }x^{-\rho_j -\frac {2 i \pi \nu}{\log q}} \omega_j 
+ d_{\mathcal{F}} B^t
$$ where $(B^t)_{t\in \R}$ is a smooth family of smooth leafwise functions on $S$ continuous on $S.$
\end{lemma}
\begin{proof} One considers only the case $\nu=0$ in 
order to simplify the notations. Using the proof of the Poincar\'e lemma and the fact that 
$\psi^0_x={\rm Id},$ 
one can write along the leaves of $\mathcal{L}:$
$$
\forall x\in \R^{+*},\; (\psi^t_x)^*(\omega_j x^{-\rho_j} )= 
\omega_j x^{-\rho_j} + d_{\mathcal{L}} A^{t,x}
$$ where $ (A^{t,x})$ is a smooth family of smooth leafwise functions on $\mathcal{L}.$
Now we use the $\beta_l$ and $I_l$  ($l=1,2$) introduced in Subsection 6.1.
For $x\in \R^{+*}$ such that $q^{-k}x \in I_l$ with $k\in \Z$, we set 
$B_l^{t,x}=  e^{ t \rho_j} \beta_l(q^{-k} x) (q^{-k})^* A^{t, q^{-k}x}$. Otherwise we set 
$B_l^{t,x}=0.$ Then, using the formula $q^* \omega_j= q^{\rho_j} \omega_j,$ one gets the result by setting 
$B^t= B_1^{t,x} + B_2^{t,x}.$
\end{proof}
Using Proposition \ref{Prop:basis} and the previous 
Lemma one then checks that

$$
{\rm TR} \int_\R \alpha(t) e^{t \Theta_1}  d t \,=\,
\sum_{j=1}^{2 g}\sum_{\nu \in \Z}   \int_\R \alpha(t) e^{t\, ( \rho_j+\frac {   2 i \pi \nu }{\log q}) }  d t.
$$

\subsection{Proof of Theorem \ref{Thm:trace}. 3] } $\;$

Notice that in the previous subsection we have shown that 
for each $j\in \{0,1,2\}$,  the operator 
$\int_{\R} \alpha(s) \pi^j_\tau \circ \phi^*_s\, ds \circ \pi^j_\tau$ is trace class.
 
 Jesus Alvarez-Lopez pointed out to us the following lemma:
\begin{lemma} (Alvarez-Lopez) \label{lem:jesus}
Set $\Theta (l,x)= \phi^{-\log x} (l,1)$ and $q (l,x)= (q\cdot l, q  x)$ where $(l,x)\in \mathcal{L} \times \R^{+*}.$ Then 
one has:

\noindent $ \Theta^{-1}\circ \phi^t \circ \Theta (l,x)=  (l, x e^{-t})$. Moreover, 
for any $(l,x)\in \mathcal{L}\times \R^{+*},$ one has 
$\Theta^{-1}\circ q \circ \Theta( l, x)= (\psi^{\log q}_q \circ q (l), q x).$
 
\end{lemma}

\medskip This Lemma shows that it is enough to prove Theorem \ref{Thm:trace}. 3] when
the flow $\phi^t$ is of the form  $\phi^t(l,x)=(l, x e^{-t})$, by abuse of notation we shall still write
$\mathcal{R}_t$, $\pi^j_\tau$ 
 instead of $\Theta^{*}\circ \mathcal{R}_t \circ( \Theta^{-1})^*$, 
$\Theta^{*}\circ \pi^j_\tau \circ( \Theta^{-1})^*$ ....etc.

\begin{proposition} \label{Prop:=}
For each $j\in \{0,1,2\}$,  denote by $\Delta^j_{\mathcal{F}}$ the leafwise signature 
Laplacian of $(S=\frac{\mathcal{L} \times \R^{+*}}{ q^\Z},\mathcal{F})$ acting on leafwise differential forms of degree $j.$
{\item 1]} For any $(t,t^\prime) \in \R^{+*} \times \R^{+*}$ and any $\alpha \in C^\infty_c(\R , \R)$,  the operator 
$$
\int_\R \alpha(s)( \phi^s)^* d s\circ \mathcal{R}_t \circ e^{-t^\prime \Delta^j_{\mathcal{F}}}
$$ is trace class.
{\item 2]} For any $(t,t^\prime) \in \R^{+*} \times \R^{+*},$ one has:
$$
\sum_{j=0}^2 (-1)^j {\rm TR}  \int_\R \alpha(s) (\phi^s)^* d s  \circ \mathcal{R}_t \circ e^{-t^\prime \Delta^j_{\mathcal{F}}}\,=\, \sum_{j=0}^2 (-1)^j {\rm TR} \int_\R \alpha(s) (\phi^s)^* d s \circ \pi^j_\tau.
$$
\end{proposition}
\begin{proof} 1] In the case of a standard Riemannian foliation, Alvarez-Lopez and Kordyukov have proved (see \cite{A-K3}  \cite{A-K2}) that 
$\int_\R \alpha(s) (\phi^s)^* d s \circ e^{-t^\prime \Delta^j_{\mathcal{F}}}
$ is trace class. Their proof can be adapted in our situation (with a 
$p-$adic transversal). 

2] Proceeding as in the proof of Lemma 3.3 of \cite{A-K2} one checks 
easily that
$$
\sum_{j=0}^2 (-1)^j {\rm TR}  \int_\R \alpha(s) (\phi^s)^* d s  \circ \mathcal{R}_t \circ e^{-t^\prime \Delta^j_{\mathcal{F}}}
$$ does not depend on $t^\prime >0$. 
Moreover, proceeding as in the proof of Lemma 3.2 of \cite{A-K2} one checks  that for fixed $t>0$
$$
\lim_{t^\prime \rightarrow +\infty} \sum_{j=0}^2 (-1)^j {\rm TR}  \int_\R \alpha(s) (\phi^s)^* d s  \circ \mathcal{R}_t \circ e^{-t^\prime \Delta^j_{\mathcal{F}}}=\sum_{j=0}^2 (-1)^j {\rm TR}  \int_\R \alpha(s) (\phi^s)^* d s  \circ \mathcal{R}_t \circ \pi^j_\tau.
$$ Since, by Proposition \ref{Prop:cte}: $\mathcal{R}_t \circ \pi^j_\tau= \pi^j_\tau$, one gets the result.
\end{proof}
\begin{proposition} Assume that the support of $\alpha$ is contained 
in $[-\frac{\log q}{2}, \frac{\log q}{2}]$, then:
$$
\sum_{j=0}^2 (-1)^j {\rm TR} \int_\R \alpha(s) (\phi^s)^* d s \circ \pi^j_\tau= \chi_\Lambda (\mathcal{F}) \alpha(0).
$$ (See the notations of Theorem \ref{Thm:trace}. 3].).
\end{proposition}
\begin{proof} Observe that the length of a closed orbit of $\phi^t$ 
is at least $\log q.$ One has just to adapt the arguments of the proof 
of Theorem 1.2 and Proposition 4.3 of \cite{A-K2}. 
\end{proof}
\begin{proposition} Assume that the support of $\alpha$ does not meet 
 $[-\frac{\log q}{4}, \frac{\log q}{4}]$, then:
$$
\sum_{j=0}^2 (-1)^j {\rm TR} \int_\R \alpha(s) (\phi^s)^* d s \circ \pi^j_\tau= 
\sum_{j=0}^2 (-1)^j  \int_S \int_\R \alpha(s) {\rm Tr}_j (D\phi^s)^* \delta_{\phi^s(y)=y} \, d \mu (y) d s
 $$ where ${\rm Tr}_j $ denotes the bundle endomorphism trace of $ (D\phi^t)^*$ acting 
on leafwise exterior forms of degree $j.$ The measure $\mu$ is defined 
in Proposition \ref{Wilson}.3]
\end{proposition}
\begin{proof} In the left hand side of the equality of Proposition 
\ref{Prop:=} one lets $t$ and $t^\prime$ go to $0^+$ and one gets 
the desired result.
\end{proof}
From the previous results we deduce that Theorem \ref{Thm:trace}.3]  follows from the:
\begin{proposition} Assume that the support of $\alpha$ does not meet 
 $[-\frac{\log q}{4}, \frac{\log q}{4}]$, then:
$$
\sum_{j=0}^2 (-1)^j \int_S \int_\R \alpha(s) {\rm Tr}_j (D\phi^s)^* \delta_{\phi^s(y)=y} \, d \mu (y) d s
=\sum_{\gamma} \sum_{k \geq 1} l(\gamma)
\left( \,e^{- k l(\gamma) } \alpha( -k l(\gamma) ) + \alpha( k l(\gamma) )\, \right)
$$ where 
$\gamma$ runs over the set of primitive closed orbits of $\phi^t$ 
and $l(\gamma)$ denotes the length of $\gamma.$

\end{proposition}
\begin{proof}   We recall Lemma \ref{lem:jesus} and the fact 
that we are reduced to the case $\phi^t(l,x)=(l,x e^{-t}).$ Assumptions (i) and (iii) of Section \ref{Section:curve} allows to see that $y \in \gamma,$ ${\rm det}\, D \phi^{l(\gamma)}(y)_{| T_y \mathcal{F}}>0.$ Moreover, using these Assumptions (i) and (iii) (especially \eqref{eq:jac})  and proceeding exactly as in the proof 
of Lemma \ref{lem:trans} one gets the proposition.
\end{proof}

\noindent One has proved the assertion $\Re \rho_j = \frac{1}{2}$ 
in Proposition \ref{Prop:basis}  as a consequence of Assumption (ii) of Section \ref{Section:curve}. 
Nevertheless, following an argument of Serre \cite{Serre} and of Deninger-Singhof [Prop 4.6]\cite{D-S} 
we are going to explain how the assertion $\Re \rho_j = \frac{1}{2}$ follows formally 
from the equality $(\phi^t)^*[\lambda_g]=e^t\,[\lambda_g ]$ of Proposition \ref{Wilson}.2].

We replace the hermitian scalar product on $H^1_\tau$ introduced before 
Theorem \ref{Thm:trace} by the following (equivalent) one.
$$
\forall \alpha_1, \alpha_2 \in H^1_\tau,\quad 
< \alpha_1 ; \alpha_2> = - (\alpha_1 \cup J \overline{ \alpha_2} ; C(\Lambda) )
$$ where the operator $J$ is multiplication by 
$\sqrt{-1}$ [resp. $-\sqrt{-1}$] on $H^{1,0}_\tau$ [resp. $H^{0,1}_\tau$].
\begin{lemma} 
{\item1]} For any $\nu \in \Z\setminus[0\},$ 
$(x^{\frac{2 i \pi \nu} {\log q}} \lambda_g ; C(\Lambda))=0.$ But 
$(\lambda_g ; C(\Lambda)) \not=0.$
{\item2]} For any $\alpha_1, \alpha_2 \in H^1_\tau$ one has:
$$
\forall t \in \R,\; < e^{t \Theta_1} \alpha_1 ; e^{t \Theta_1} \alpha_2 > =
( (\phi^t)^*(\alpha_1 \cup J \overline{ \alpha_2} ) ; C(\Lambda) ).
$$
{\item3]} For any $\alpha_1, \alpha_2 \in H^1_\tau$ one has:
$$
\forall t \in \R,\; < e^{t \Theta_1} \alpha_1 ; e^{t \Theta_1} \alpha_2 > =
e^t < \alpha_1 ; \alpha_2>.
$$
\end{lemma}
\begin{proof} {\item 1]} Easy computation.
{\item 2]} Recall that we have defined two notions of leafwise reduced cohomology 
(see Definition \ref{def:form} and the proof of Theorem \ref{Theorem:Hodge}).  We have  a natural map between them:
$$
\overline{H}^j_\mathcal{F}(S)  \rightarrow \overline{H}^j_{0, +\infty}.
$$
 Observe that $e^{t\Theta_1} \alpha= (\phi^t)^* \alpha$ modulo 
$\overline{{\rm Im}\, d_{\mathcal{F}}}$ and that the Ruelle-Sullivan current 
is closed. One then checks that 
$$
< e^{t \Theta_1} \alpha_1 ; e^{t \Theta_1} \alpha_2 > =
( (\phi^t)^* \alpha_1  \cup J\, \overline{ (\phi^t)^* \alpha_2}  ; C(\Lambda) ).
$$
Since $(\phi^t)^*$ commutes with $J$ and the complex conjugation, one gets the result.
{\item 3]} One can write 
$$
 \pi^2_\tau (\alpha_1 \cup J \overline{ \alpha_2})= \sum_{\nu \in \Z} c_\nu 
x^{\frac{2 i \pi \nu} {\log q}} \lambda_g \in H^2_\tau.
$$ 
 Then, using parts 1] and  2] one checks that
$ < e^{t \Theta_1} \alpha_1 ; e^{t \Theta_1} \alpha_2 > = 
((\phi^t)^*(c_0 \lambda_g)  ; C(\Lambda) )$ and

 $<  \alpha_1 ;  \alpha_2 > = ( c_0 \lambda_g  ; C(\Lambda) )$. Since $(\phi^t)^*[\lambda_g]= e^t [\lambda_g] $ by Proposition 
\ref{Wilson}.2] one gets the result.
\end{proof}

Now part 3] of the previous Lemma implies that
$$
\frac{d }{ d t} \left(< e^{t \Theta_1} \alpha_1 ; e^{t \Theta_1} \alpha_1 >\right)_{| t=0}=
< \alpha_1 ; \alpha_1>.
$$
Therefore one has:
$$
< (\Theta_1 -\frac{1}{2}) (\alpha_1) ; \alpha_1> + <\alpha_1 ; (\Theta_1 -\frac{1}{2}) (\alpha_1)> = 0.
$$ Now if $\alpha_1 \in H^1_\tau \setminus\{0\}$ is such that 
$\Theta_1(\alpha_1) = \rho \,\alpha_1$ with $\rho \in \C$ 
then one gets $\rho-\frac{1}{2}=-(\overline{\rho}-\frac{1}{2}).$
 Thus one gets $\Re\, \rho = \frac{1}{2}$ which from the equality 
$(\phi^t)^*[\lambda_g]=e^t [\lambda_g]$ as desired.
\smallskip
\section{Appendix: Renormalization group flow}$\;$

We first briefly and informally recall Wilson's view point following 
[pages 554 and 557]\cite{Del}. 
We consider a set $\mathcal{S} $ of QFT defined by lagrangians. Let $0 <\Lambda_0 < \Lambda_1$ 
be two scales of momenta. For each theory $L\in \mathcal{S}$, one finds another theory 
$R_{\Lambda_1, \Lambda_0}(L)$ which is the effective theory at the scale 
$\Lambda_1$ of the original theory $L$ at the scale $\Lambda_0.$

In terms of Feynman integrals for QFT defined on $\R^n$ one can write:
$$
\int_{\mathcal{B}(\Lambda_0)} A(\phi) \left( \int_{\mathcal{C}(\Lambda_0, \Lambda_1)} e^{-L(\phi + \eta)} D \eta \right) D \phi=
\int_{\mathcal{B}(\Lambda_0)} A(\phi) e^{-R_{\Lambda_1, \Lambda_0}(L)\,(\phi)} D \phi
$$ where $A$ is a function of the field $\phi$, and 
$ \mathcal{C}(\Lambda_0, \Lambda_1)$ (resp. $\mathcal{B}(\Lambda_0) $) denotes the set of fields whose Fourier transform has support in the corona 
$\{\xi \in \R^n,\; \Lambda_0 \leq | \xi| \leq \Lambda_1 \}$ (resp. the ball $ \{ \xi \in \R^n,\;  |\xi | \leq \Lambda_0\}$).

Then the renormalization (semi-)group flow is defined by:
$$
\forall t \in \R^+,\; \forall ( L, \Lambda_1) \in  \mathcal{S}\times\R^{+*} ,\; \phi^t(L, \Lambda_1)=
( \,
 R_{\Lambda_1, e^{-t} \Lambda_1}(L), \Lambda_1 e^{-t} \,).
$$ Notice that if $L$ is a free lagrangian then  
 $\phi^t( L, \Lambda_1)=(  L, \Lambda_1 e^{-t})$ for any $t\geq 0.$

\smallskip
In \cite{CK0} and \cite{CK1} Connes and Kreimer have developed a mathematical 
theory of renormalization of perturbative QFT. Let $G$ be 
the group of characters  of the dual Hopf algebra 
of the enveloping algebra associated with the $1-$particle irreducible Feynman graphs. The unrenormalized theory gives rise to a 
meromorphic loop $\gamma(z) \in G, \,z \in \C.$ Connes and Kreimer 
have shown that the renormalized theory is the evaluation 
at the integer dimension $z_0$ of space-time of the holomorphic 
part $\gamma_+$ of the Birkhoff decomposition of $\gamma$ for the Riemann-Hilbert problem.
Moreover they view the renormalization group 
as a subgroup of  $G$. Then, for massless  QFT  they recover the action of the 
renormalization group on lagrangians. 
 
\medskip
\section{Acknowledgements}$\;$

I thank  C. Deninger whose helpful comments 
on the first version have allowed to improve  this paper.
I thank   A. Abdesselam, J. Alvarez-Lopez, A. Besser, U. Buenke,  L. Boutet de Monvel,  C. Deninger, D. Harari, S. Haran, M. Marcolli,  G. Skandalis and, D. Sullivan 
for fruitful   discussions.
Part of this work was done while the author was visiting the universities 
of Bar Ilan, Beer Sheva, Technion, Tel Aviv and Savoie: he would like to thank these 
institutions for their very warm hospitality.
I thank Etienne Blanchard for his interest in this work and having given me the opportunity to present these 
results in the $C^*-$algebras seminar of Paris.

	\bibliographystyle{amsalpha}

\begin{thebibliography}{ACDE}

\bibitem[A-K00]{A-K2} J. A. Alvarez Lopez and Y. Kordyukov: {\it Distributional 
Betti numbers of transitive foliations of codimension one}. Foliations: geometry 
and dynamics (Warsaw, 2000), World Sci. Publishing, River Edge, NJ, (2002), 
pages 159-183.

\bibitem[A-K01]{A-K1} J. A. Alvarez Lopez and Y. Kordyukov: {\it Long time 
behaviour of leafwise heat flow for Riemannian foliations}. Compositio Math. 
{\bf 125} (2001), pages 129-153.

\bibitem[A-K05]{A-K3} J. A. Alvarez Lopez and Y. Kordyukov: {\it
Distributional Betti numbers for Lie foliations}. Preprint.

\bibitem[Co94]{C 1} A. Connes: {\it \underline{Noncommutative Geometry}}. 
Academic Press.

\bibitem[Co99]{C 2} A. Connes: {\it Trace formula in noncommutative
geometry and the zeroes of the Riemann zeta function}. Selecta Math.
(N.S.) 5 (1999), no. 1, 29-106.

\bibitem[Co00] {C3} A. Connes: {\it Noncommutative Geometry Year 2000}, 
Special Volume GAFA 2000 Part II, pages 481-559.

\bibitem[Co02]{C4} A. Connes:{\it Sym\'etries Galoisiennes et Renormalisation}, 
S\'eminaire Bourbaphy,  Octobre 2002, pages 75-91.



\bibitem[CK00]  {CK0} A. Connes, and D. Kreimer:{\it Renormalization in quantum field 
theory and the Riemann Hilbert problem I}, Comm. Math. Phys. {\bf 210}, 2000, No 1, pages 
249-273.

\bibitem[CK01]  {CK1} A. Connes, and D. Kreimer:{\it Renormalization in quantum field 
theory and the Riemann Hilbert problem II, the $ \beta$ function, 
diffeomorphisms and the renormalization group}, Comm. Math. Phys. {\bf 216},2001, No 1, pages 
249-273.

\bibitem[Co-Ma04a] {CM1} A. Connes, and M. Marcolli:{\it 
$\Q-$lattices: quantum statistical mechanics and Galois theory}, to appear 
in Journal of Geometry and Physics.

\bibitem[Co-Ma04b] {CM2} A. Connes, and M. Marcolli:{\it From Physics to Number theory via Noncommutative Geometry. Part II:  Renormalization, the Riemann-Hilbert correspondence, and 
motivic Galois theory}, to appear in the volume "Frontiers in Number Theory, Physics, and Geometry".


\bibitem[QFT]{Del} P. Deligne and Al: {\it \underline{Quantum Fields and Strings:
A Course for Mathematicians, vol I}}. Amer. Math. Soc.  Institute for Advanced Study (1996).




\bibitem[De94]{D 3} C. Deninger: {\it Motivic $L-$functions 
and regularized determinants }. Proc. Symp. Pure Math. {\bf 55}, 1 (1994), pages 707-743. 


\bibitem[De98]{D 1} C. Deninger: {\it Some analogies between number
theory and dynamical systems onf foliated spaces}.
Doc. Math. J. DMV. Extra volume ICM I, (1998), pages 23-46.

\bibitem[De99]{D 0c} C. Deninger: {\it On dynamical systems and their possible significance 
for Arithmetic Geometry}. In: A. Reznikov, N. Schappacher (eds.), Regulators in Analysis, 
Geometry and Number Theory. Progress in Mathematics {\bf 171}, (1999), 
Birkhauser, pages 29-87. 


\bibitem[De-Si00] {D 1b}C. Deninger and W. Singhof: {\it A note on dynamical trace formulas}, 
Dynamical, spectral and arithmetic zeta functions (San Antonio, TX, 1999), Contemp. 
Math., 290, Amer. Math. Soc., Providence, RI, (2001), pages 41-55.

\bibitem[De-Si02] {D-S} C. Deninger and W. Singhof: {\it Real polarizable Hodge structures 
arising from foliations}.  Annals of Global Analysis and Geometry 21,  2002,  pages377-399

\bibitem[De01]{D 2b} C. Deninger: {\it Number theory and dynamical systems 
on foliated spaces}. Jahresberichte der DMV. {\bf 103},  (2001), No 3, pages 
79-100.

\bibitem[De01b]{D 1c} C. Deninger: {\it A note on arithmetic topology and dynamical systems}. 
Algebraic number theory and algebraic geometry, Contemp. Math., 300, Amer. Math. Soc., 
Providence, RI, (2002), pages 99-114.



\bibitem[De02]{D 2} C. Deninger: {\it On the nature of explicit
formulas in analytic number theory, a simple example}. Number theoretic methods (Iizuka, 2001),
 Dev. Math., 8, Kluwer Acad. Publ., Dordrecht, (2002), pages 97-118.

\bibitem[Der04]{Deroin} B. Deroin:{\it Non rigidity of hyperbolic Riemann surfaces 
lamination}. Preprint 2004.
 
\bibitem[Ghys99]{Ghys} E. Ghys:{\it Laminations par surfaces de Riemann}.
 Dynamique et g\'eom\'etrie complexes (Lyon, 1997),  
Panor. Synthses, 8, Soc. Math. France, Paris, (1999), pages 49-95.
\bibitem[G-S77]{G-S} V. Guillemin and S. Sternberg:\underline{{\it Geometric 
Asymptotics}}. Mathematical surveys and monographs. Number 14.
Published by the A.M.S. 

\bibitem[Har05]{Haran} S. Haran:{\it Arithmetic as Geometry I. The language of 
non-additive geometry}. Preprint  2005.

\bibitem[Ih73]{Ihara} Y. Ihara:{\it On $(\infty\times p)-$adic coverings 
of curves (the simplest example)}. Trudy Mat. Inst. Steklov {\bf 132} (1973),
 pages 133-148.

\bibitem[Illes99]{Illies} G. Illies:{\it Cramer Functions and Guinand Equations}. 
Preprint IHES (1999).



\bibitem[Lei03] {El} E. Leichtnam: {\it An invitation to 
Deninger's work on arithmetic zeta functions}.  Geometry, spectral theory, groups, and dynamics, Contemp. Math., 387, Amer. Math. Soc., Providence, RI, (2005), pages 
201-236.

\bibitem[Lei06] {El1} E. Leichtnam: {\it Renormalization group flow 
and arithmetic Zeta functions}. In preparation.


\bibitem[Meyer03]{Meyer} R. Meyer:{\it On a representation of the Idele class group related 
to primes and zeros of $L-$functions}, Preprint 2003.
\bibitem[Mil80] {Milne} J. Milne:{\it \underline{ \'Etale Cohomology}}, Princeton 
Mathematical Series, 33, (1980).

\bibitem[Oor73]{Oort} F. Oort: {\it Lifting an endomorphism of an elliptic curve to 
characteristic zero}. Indag. Math. 35, (1973), pages 466-470.

\bibitem[Po84] {Polchinski} J. Polchinski: {\it Renormalization and effective lagrangians}. 
Nuclear Physics B231, (1984), pages 269-295. 




\bibitem[Se]{Serre0} J-P. Serre: {\it \underline{Groupes alg\'ebriques et corps de classe}}. Editions
Hermann (1959).
\bibitem[Se60]{Serre} J-P. Serre: {\it Analogues Kaehleriens de certaines 
conjectures de Weil}. Annals of Math. 65, (1960), pages 392-394.
\bibitem[Si92] {Silverman} J. Silverman: {\it \underline{The arithmetic of elliptic curves}}. 
Graduate text in Math. 106  (1992).
\bibitem[Sul93]{Sullivan} D. Sullivan:{\it Linking the universalities of Milnor-Thurston, Feigenbaum and Ahlfors-Bers}.  Topological methods in modern mathematics (Stony Brook, NY, 1991),   Publish or Perish, Houston, TX, (1993), pages 543-564.

\end{thebibliography}

\end{document}